\documentclass[11pt]{amsart} \textwidth=14.5cm \oddsidemargin=1cm
\usepackage[rgb]{xcolor}
\usepackage{tikz}
\usepackage{verbatim}
\usepackage{amsmath}
\usepackage{amsxtra}
\usepackage{amscd}
\usepackage{amsthm}
\usepackage{amsfonts}
\usepackage{amssymb}
\usepackage{eucal}
\usepackage{float}
\usepackage{mathrsfs}
\usepackage{graphicx}
\usepackage{stmaryrd}
\usepackage{tikz-cd}
\usetikzlibrary{arrows}
\usetikzlibrary{cd}
\usepackage{ytableau}

\usepackage{latexsym,todonotes}

\usepackage[linktocpage=true]{hyperref}
\hypersetup{colorlinks,linkcolor=blue,urlcolor=orange,citecolor=blue}
\usepackage[all, cmtip]{xy}

\usepackage{stmaryrd}
\usepackage{color} 
\usepackage{amsrefs}

\newtheorem{thm}{Theorem} [section]

\newtheorem{lem}[thm]{Lemma}
\newtheorem{prop}[thm]{Proposition}

\theoremstyle{definition}
\newtheorem{definition}[thm]{Definition}

\theoremstyle{remark}
\newtheorem{rem}[thm]{Remark}

\numberwithin{equation}{section}
\begin{document}

%Referring commands:
\newcommand{\thmref}[1]{Theorem~\ref{#1}}
\newcommand{\secref}[1]{Section~\ref{#1}}
\newcommand{\lemref}[1]{Lemma~\ref{#1}}
\newcommand{\propref}[1]{Proposition~\ref{#1}}
\newcommand{\corref}[1]{Corollary~\ref{#1}}
\newcommand{\remref}[1]{Remark~\ref{#1}}
\newcommand{\eqnref}[1]{(\ref{#1})}

\newcommand{\exref}[1]{Example~\ref{#1}}

%Theorem for the introduction
\newtheorem{innercustomthm}{{\bf Theorem}}
\newenvironment{customthm}[1]
  {\renewcommand\theinnercustomthm{#1}\innercustomthm}
  {\endinnercustomthm}
  
  \newtheorem{innercustomcor}{{\bf Corollary}}
\newenvironment{customcor}[1]
  {\renewcommand\theinnercustomcor{#1}\innercustomcor}
  {\endinnercustomthm}
  
  \newtheorem{innercustomprop}{{\bf Proposition}}
\newenvironment{customprop}[1]
  {\renewcommand\theinnercustomprop{#1}\innercustomprop}
  {\endinnercustomthm}

\newcommand{\bbinom}[2]{\begin{bmatrix}#1 \\ #2\end{bmatrix}}
\newcommand{\cbinom}[2]{\set{\^!\^!\^!\begin{array}{c} #1 \\ #2\end{array}\^!\^!\^!}}
\newcommand{\abinom}[2]{\ang{\^!\^!\^!\begin{array}{c} #1 \\ #2\end{array}\^!\^!\^!}}
\newcommand{\qfact}[1]{[#1]^^!}

%Simplified symbols:
\newcommand{\nc}{\newcommand}

\nc{\Ord}{\text{Ord}_v}

 \nc{\A}{\mathcal A} 
  \nc{\G}{\mathbb G} 
\nc{\Ainv}{\A^{\rm inv}}
\nc{\aA}{{}_\A}
\nc{\aAp}{{}_\A'}
\nc{\aff}{{}_\A\f}
\nc{\aL}{{}_\A L}
\nc{\aM}{{}_\A M}
\nc{\Bin}{B_i^{(n)}}
\nc{\dL}{{}^\omega L}
\nc{\Z}{{\mathbb Z}}
 \nc{\C}{{\mathbb C}}
 \nc{\N}{{\mathbb N}}
 \nc{\R}{{\mathbb R}}
  \renewcommand{\S}{{\mathcal S}}
  \nc{\I}{{\mathcal I}}
 \nc{\fZ}{{\mf Z}}
 \nc{\F}{{\mf F}}
 \nc{\Q}{\mathbb{Q}}
 \nc{\la}{\lambda}
 \nc{\ep}{\epsilon}
 \nc{\h}{\mathfrak h}
 \nc{\He}{\bold{H}}
 \nc{\htt}{\text{tr }}
 \nc{\n}{\mf n}
 \nc{\g}{{\mathfrak g}}
 \nc{\DG}{\widetilde{\mathfrak g}}
 \nc{\SG}{\breve{\mathfrak g}}
 \nc{\is}{{\mathbf i}}
 \nc{\V}{\mf V}
 \nc{\bi}{\bibitem}
 \nc{\E}{\mc E}
 \nc{\ba}{\tilde{\pa}}
 \nc{\half}{\frac{1}{2}}
 \nc{\hgt}{\text{ht}}
 \nc{\ka}{\kappa}
 \nc{\mc}{\mathcal}
 \nc{\mf}{\mathfrak} 
 \nc{\hf}{\frac{1}{2}}
\nc{\ov}{\overline}
\nc{\ul}{\underline}

\nc{\xx}{{\mf x}}
\nc{\id}{\text{id}}
\nc{\one}{\bold{1}}
\nc{\mfsl}{\mf{sl}}
\nc{\mfgl}{\mf{gl}}
\nc{\ti}[1]{\textit{#1}}
\nc{\Hom}{\text{Hom}}
\nc{\Cat}{\mathscr{C}}
\nc{\CatO}{\mathscr{O}}
\renewcommand{\O}{\mathscr{O}}
\nc{\Tan}{\mathscr{T}}
\nc{\Umod}{\mathscr{U}}
\nc{\Func}{\mathscr{F}}
\nc{\Kh}{\text{Kh}}
\nc{\Khb}[1]{\llbracket #1 \rrbracket}

\nc{\ua}{\mf{u}}
\nc{\nb}{u}
\nc{\inv}{\theta}
\nc{\mA}{\mathcal{A}}
\newcommand{\TT}{\mathbf T}
\newcommand{\TA}{{}_\A{\TT}}
%\nc{\bV}{{\bold V}}
\newcommand{\tK}{\widetilde{K}}
\newcommand{\al}{\alpha}
\newcommand{\Fr}{\bold{Fr}}

%quantum groups U
\nc{\Qq}{\Q(v)}
\nc{\U}{\bold{U}}
\nc{\uu}{\mathfrak{u}}
\nc{\Udot}{\dot{\U}}
\nc{\f}{\bold{f}}
\nc{\fprime}{\bold{'f}}
\nc{\B}{\bold{B}}
\nc{\Bdot}{\dot{\B}}
%\nc{\psi}{\psi}
\nc{\Dupsilon}{\Upsilon^{\vartriangle}}
\newcommand{\T}{\texttt T}
\newcommand{\vs}{\varsigma}
\newcommand{\Pa}{{\bf{P}}}
\newcommand{\Padot}{\dot{\bf{P}}}

%Ui
\nc{\ipsi}{\psi_{\imath}}
\nc{\Ui}{{\bold{U}^{\imath}}}
\nc{\uidot}{\dot{\mathfrak{u}}^{\imath}}
\nc{\Uidot}{\dot{\bold{U}}^{\imath}}
%\nc{\aA \Ui}{{{}_\A\bold{U}^{\imath}}}
%\nc{\aA \Uidot}{{}_\A\dot{\bold{U}}^{\imath}}
 \nc{\be}{e}
 \nc{\bff}{f}
 \nc{\bk}{k}
 \nc{\bt}{t}
 \nc{\bs}{\backslash}
 \nc{\BLambda}{{\Lambda_{\inv}}}
\nc{\Ktilde}{\widetilde{K}}
\nc{\bktilde}{\widetilde{k}}
\nc{\Yi}{Y^{w_0}}
\nc{\bunlambda}{\Lambda^\imath}
\newcommand{\Iwhite}{\I_{\circ}}
\nc{\ile}{\le_\imath}
\nc{\il}{<_{\imath}}

%newUi
\newcommand{\ff}{B}

%bilinear form
%\nc{\qq}{(q_i^{-1}-q_i)}
%\nc{\qqq}{(1-q_i^{-2})^{-1}}
%\nc{\qqqj}{(1-q_j^{-2})^{-1}}

%Blackdots
\nc{\etab}{\eta^{\bullet}}
\newcommand{\Iblack}{\I_{\bullet}}
\newcommand{\wb}{w_\bullet}
\newcommand{\UIblack}{\U_{\Iblack}}

%color
\newcommand{\blue}[1]{{\color{blue}#1}}
\newcommand{\red}[1]{{\color{red}#1}}
\newcommand{\green}[1]{{\color{green}#1}}
\newcommand{\white}[1]{{\color{white}#1}}

%i-divided powers
%\newcommand{\dvev}[1]{{\mathfrak{t}}_{\ev}^{{(#1)}}}
%\newcommand{\dv}[1]{{\mathfrak{t}}_{\odd}^{{(#1)}}}
\newcommand{\dvd}[1]{t_{\odd}^{{(#1)}}}
\newcommand{\dvp}[1]{t_{\ev}^{{(#1)}}}
\newcommand{\ev}{\mathrm{ev}}
\newcommand{\odd}{\mathrm{odd}}

\newcommand\TikCircle[1][2.5]{{\mathop{\tikz[baseline=-#1]{\draw[thick](0,0)circle[radius=#1mm];}}}}

\newcommand{\commentcustom}[1]{}

\raggedbottom

\title[Transfer using Fourier transform and minimal representation of $E_7$]
{Transfer using Fourier transform and minimal representation of $E_7$}

\author[Nhat Hoang Le and Bryan Peng Jun Wang]{Nhat Hoang Le and Bryan Peng Jun Wang}
\address{Department of Mathematics, Block S17, National University of Singapore, 10 Lower Kent Ridge Drive, 119076.}
\email{lnhoang@nus.edu.sg}
\address{Department of Mathematics, Harvard University, Cambridge, MA 02138, USA}
\email{bwangpengjun@math.harvard.edu}

\begin{abstract}
    In this paper, we study the Sakellaridis-Venkatesh conjecture for the rank-1 spherical variety $X=\text{Spin}_9\backslash F_4$ using an exceptional theta correspondence. We establish the correct transfer map satisfying relative character identities in this case and show that our transfer map agrees with the formula in \cite{S18}. We also formulate the local relative characters for the degenerate Whittaker period of \cite{MWZ26} associated with $X$. Moreover, we show how our techniques lead to a characterization of $X$-relatively cuspidal representations.
\end{abstract}

%\vspace{.3cm}
\maketitle

\setcounter{tocdepth}{1}
\tableofcontents

\section{Introduction}

\subsection{The relative Langlands program}

The theory of periods of automorphic forms has been an important topic of study in the Langlands program. In particular, the nonvanishing of periods of automorphic forms on a reductive group $G$ relative to a certain subgroup $H$ of $G$ can be related to special values of certain automorphic $L$-functions. In the local field setting, the corresponding problem is to give a characterisation of irreducible smooth representations of $G$ whose spaces of $H$-invariant functionals are nontrivial. By Frobenius reciprocity, this is equivalent to classifying $G$-submodules of $C^\infty(X)$, where $X=H\backslash G$. We call such representations $X$-distinguished. We can also consider an $L^2$-analog of this question, namely, to give a spectral decomposition of the unitary representation $L^2(X)$.

In \cite{SV17}, Sakellaridis and Venkatesh give precise predictions for these problems in the context of spherical varieties, which is now regarded as the relative Langlands program. To be more precise, when $X=H\backslash G$ is a spherical variety, following the spirit of the Langlands philosophy, they associate the following data:
\begin{itemize}
    \item a Langlands dual group $X^\vee$ (assuming $G$ is split) and a canonical map (up to conjugacy):
    $$
    \iota: X^\vee\times \text{SL}_2(\mathbb{C}) \rightarrow G^\vee.
    $$
    \item a graded finite-dimensional representation $V_X$ of $X^\vee$ (which is expected to be symplectic), which gives rise to an $L$-function $L_X(s,\rho)$ for each $L$-parameter $\rho$ valued in $X^\vee$.
\end{itemize}
They conjectured that $X$-distinguished representations of $G$ belong to $A$-packets whose associated $A$-parameters factor through $\iota$. In other words, they are Langlands functorial lifted via $\iota$ from a split group $G_X$ whose dual is $X^\vee$. Often, it is helpful to think of the $X$-distinguished representations of $G$ as being lifted from the Whittaker variety $(N_X,\psi)\backslash G_X$ -- in other words, from generic representations of $G_X$ -- rather than the group $G_X$ itself. In the global setting, the relevant period integral is controlled by the automorphic L-function $L_X$. More recently, by considering the cotangent space $T^*X$ instead of the spherical variety $X$, Ben-Zvi, Sakellaridis and Venkatesh \cite{BZSV24} extended this program to the setting of hyperspherical varieties. However, since the framework in \cite{SV17} has better fit to the case considered in this paper, we will restrict our focus to the setting of spherical varieties.

We now give a more precise discussion of the Sakellaridis-Venkatesh conjectures over local fields (including both smooth and $L^2$-settings) and global fields, under some simplifications which are convenient in the case of our interest.

In the local, smooth setting, one would like to address the question of determining $\text{Hom}_H(\pi,\mathbb{C})$, for any irreducible smooth representation $\pi$ of $G(F)$. We want to give a map 
$$
\iota_*: \text{Irr}(G_X(F)) \rightarrow \text{Irr}(G(F)),
$$
such that for any $\pi\in \text{Irr}(G(F))$, there is an isomorphism 
$$
f: \bigoplus_{\sigma:\iota_*(\sigma)=\pi}\text{Hom}_H(\sigma,\psi)\cong \text{Hom}_H(\pi,\mathbb{C}).
$$
In this setting, the Sakellaridis-Venkatesh conjecture gives a precise statement of the expectation that $X$-distinguished representations are those lifted from $G_X$. If $\iota_*$ is injective, then by the uniqueness of Whittaker models, all these Hom spaces are at most one-dimensional. This would be the case we consider in this paper. In such cases, if $\ell \in \text{Hom}_{N_X}(\sigma,\psi)$, one can define relative characters $\mathcal{B}_{\sigma,\ell}$ and $\mathcal{B}_{\iota_*(\sigma),f(\ell)}$ which are distributions on $(N_X,\psi)\backslash G_X$ and $X$ respectively. One expects to have a relative character identity relating the two distributions $\mathcal{B}_{\sigma,\ell}$ and $\mathcal{B}_{\iota_*(\sigma),f(\ell)}$.

In the $L^2$-setting, we want to give an explicit spectral decomposition of $L^2(X)$. The Sakellaridis-Venkatesh conjecture explicates the following abstract direct integral decomposition given by functional analysis
$$
L^2(X)= \int_{\Omega} \pi_\omega^{\oplus m(\omega)}d\mu_X(\omega).
$$
Namely, they conjecture that there is a map $\iota_*: \hat{G}_X \rightarrow \hat{G}$ associated to the functorial lifting $\iota$, such that we have the following unitary isomorphism
$$
L^2(X) \cong \int_{\hat{G}_X}\iota_*(\sigma)^{m(\sigma)}d\mu_{G_X}(\sigma),
$$
where $d\mu_{G_X}$ is the Plancherel measure of $G_X$ and $m(\sigma)$ is a multiplicity space which is expected to be isomorphic to the dual space of $\text{Hom}_{N_X}(\sigma,\psi)$. The above spectral decomposition gives a canonical element $\ell_{\iota_*(\sigma)}$ in $\text{Hom}_H(\iota_*(\sigma),\mathbb{C})$.

In the global setting, given a global field $k$ whose adele ring is $\mathbb{A}$, we consider the global period integral $\mathcal{P}_H:\mathcal{A}_{\text{cusp}}(G)\rightarrow \mathbb{C}$ along $H$ given by 
$$
\mathcal{P}_H(\phi)=\int_{H(F)\backslash H(\mathbb{A})}\phi(h)dh,
$$
where $\mathcal{A}_{\text{cusp}}(G)$ is the space of cusp forms on $G$. The restriction of $\mathcal{P}_H$ to a cuspidal automorphic representation $\pi=\otimes_v \pi_v$ of $G$ gives rise to an element $\mathcal{P}_{H,\pi}\in \text{Hom}_{H(\mathbb{A})}(\pi,\mathbb{C})$. One would like to address the following two questions:
\begin{enumerate}
    \item characterising those cuspidal automorphic representations $\pi$ such that $\mathcal{P}_{H,\pi}$ is nonzero via the functorial lift from $G_X$;
    \item determining whether $\mathcal{P}_{H,\pi}$ can be factorised as an Euler product of local functionals.
\end{enumerate}
Since the local Hom spaces $\text{Hom}_{H(F_v)}(\pi_v,\mathbb{C})$ are at most one dimensional for all places $v$ in the case considered in this paper, such a factorisation certainly exists. Moreover, as mentioned in the $L^2$-setting, we have a canonical basis element $\ell_{\pi_v}\in \text{Hom}_{H(F_v)}(\pi_v,\mathbb{C})$. When $\pi_v$ is unramified, Sakellaridis evaluated $\ell_{\pi_v}(\phi_{0,v})$, where $\phi_{0,v}$ is a spherical unit vector in $\pi_v$ \cite{S08,S13}. Namely, given $\pi_v \cong \iota_*(\sigma_v)$ for tempered $\sigma_v\in \hat{G}_{X,v}$, we have
$$
|\ell_{\pi_v}(\phi_{0,v})|^2=L_X^{\#}(1/2,\pi_v):=\Delta_{X,v}(0)\cdot \frac{L_{X,v}(1/2,\sigma_v)}{L(1,\sigma_v,Ad)}>0,
$$
where $\Delta_{X,v}(0)$ is a product of $L$-factors depending only on $X$ and not on the representation $\sigma_v$. The earlier-mentioned representation $V_X$ of $X^\vee$ determines what $L$-factors appear here. To obtain the Euler product factorisation of $\mathcal{P}_{H,\pi}$, we define a normalisation of $\ell_{\pi_v}$
$$
\ell_{\pi_v}^\flat=|L_X^{\#}(1/2,\pi_v)|^{-1/2}\cdot\ell_{\pi_v}.
$$
Then it remains to determine the constant $c(\pi)$ such that
$$
|\mathcal{P}_{H,v}(\phi)|^2=c(\pi)\cdot L_X^{\#}(1/2,\pi)\cdot \prod_v |\ell_{\pi_v}^\flat(\phi_v)|^2,
$$
for $\phi=\otimes_v\phi_v\in \otimes_v \pi_v$. 

In a series of works (cf. \cite{S13a,S19a,S19b,S18,S22a,S22b}), Sakellaridis studied these above problems in the context of rank-1 spherical varieties $X$. There is a classification of such spherical varieties, where rank-1 here means that the group $G_X$ is $\text{SL}_2$ or its variants, $\text{Mp}_2$ or $\text{PGL}_2$. He developed a theory of transfer of test functions from $X$ to the Whittaker variety $(N_X,\psi)\backslash G_X$ and sketched a framework to obtain the local relative character identity and a global comparison of the relative trace formula of $X$ and the Kuznetsov trace formula for $G_X$. This was motivated by considering an analogous transfer for the boundary degenerations of $X$ and $(N_X,\psi)\backslash G_X$. Moreover, he also confirmed that the given transfer map satisfies a local relative character identity in the basic cases $A_1,\,D_2$ and stated that it could be shown in the general $A_n$-case by unfolding.

On the other hand, the $L^2$-setting was treated in \cite{GG14} for essentially general rank 1 spherical varieties, whereas the smooth setting was considered in \cite{Gan19} by means of the theta correspondence. Moreover, also by using the theta correspondence, in \cite{GW18}, Gan and Wan resolved (and recovered the results of Sakellaridis) in the case $X=\text{SO}_{n-1}\backslash \text{SO}_n$.

In this paper, motivated by \cite{GW18}, we treat one of the outstanding rank-1 cases for exceptional groups -- that is the case when $X=\text{Spin}_9\backslash F_4$ -- via an exceptional theta correspondence. In this case, we have $G_X=\text{PGL}_2$. By using the exceptional theta correspondence for the dual pair $(\text{PGL}_2,F_4)$ inside $E_7$, which was recently determined completely by Karasiewicz and Savin in \cite{KS23}, we are able to give a conceptual, geometric definition of the transfer and relevant spaces of test functions. Similar to \cite{GW18}, we then establish the desired relative character identities mentioned in the smooth setting, without doing a geometric comparison. Moreover, by showing that our formula essentially agrees with one in \cite[Theorem 1.3]{S18} (Theorem \ref{thm3.5}), we confirm Sakellaridis' claim that the given operator is the correct operator of functoriality between the relative trace formula for $X$ and the Kuznetsov formula in this case. Globally, we obtain an Euler factorisation of the global period integral (Theorem \ref{thm6.6}) and the global relative character identity (Theorem \ref{thm6.7}).

\subsection{Relative cuspidality}

As an application of the local relative character identity, we discuss the notion of relative cuspidality introduced by Kato and Takano in \cite{KT08}. Namely, an $X$-distinguished representation $(\pi,V)$ of $G(F)$ is said to be $X$-relatively cuspidal if for any $v\in V$ and $\lambda\in \text{Hom}_H(\pi,\mathbb{C})$, the corresponding relative matrix coefficient
$$
F_{v,\lambda}:g\in H(F)\backslash G(F) \mapsto \lambda(\pi(g)v)
$$
is compactly supported. In the group case, when $H\subset G=H\times H$ diagonally, an H-distinguished rep is of the form $V\otimes V^*$ and $V\otimes V^*$ is $X$-relatively cuspidal if and only if $V$ is supercuspidal in the usual sense. For a general symmetric space $X$, as in \cite{Mur18}, one may ask if there exists $X$-relatively cuspidal representations $\pi$ which are not supercuspidal. There are some difficulties in addressing this question, as it is not easy to write down an explicit formula for relative matrix coefficients, let alone to check that they are compactly supported. Moreover, even when using the criterion based on Jacquet modules (see \cite[Theorem 6.2]{KT08}), it is not clear if there is a uniform way to determine the map $r_P(\lambda)$ in general. In the case when $X=\text{Spin}_9\backslash F_4$, we have an explicit isomorphism
$$
\text{Hom}_{N}(\sigma,\psi) \longrightarrow \text{Hom}_{\text{Spin}_9}(\theta_\psi(\sigma),\mathbb{C}),
$$
which yields a formula for relative matrix coefficients of $\text{Spin}_9\backslash F_4$ in terms of Whittaker coefficients of $\sigma$. This idea has been discussed in the remark at the end of section 8 in \cite{KS23}. In Theorem \ref{thm4.3}, we elaborate their discussion further by showing that $\theta_\psi(\sigma)$ is $X$-relatively cuspidal if and only if $\sigma$ is supercuspidal. This answers the above question in our case, noting that $\theta_\psi(\sigma)$ is not supercuspidal for any $\sigma \in \text{Irr}(\text{PGL}_2)$ (cf. \cite[Theorem 1.2]{KS23}).

\subsection{Degenerate Whittaker periods} In \cite{MWZ26}, the three authors proposed a relative trace formula comparison to address the numerical conjecture in \cite{BZSV24}, restated as \cite[Conjecture 1.2]{MWZ26}. More precisely, they formulated a conjecture on the degenerate Whittaker period (\cite[Conjecture 1.6]{MWZ26}) in order to reduce \cite[Conjecture 1.2]{MWZ26} to the strongly tempered case.

As noted in \cite[Remark 1.7]{MWZ26}, a key obstruction to formulating an explicit statement for \cite[Conjecture 1.6]{MWZ26} is the construction of local relative characters for these degenerate Whittaker periods. We show that the definition of local relative characters for the degenerate Whittaker period corresponding to $F_4/\text{Spin}_9$ arises naturally in the context of the exceptional theta correspondence.

\subsection{Organization of this paper}

The paper is organized as follows. In Section \ref{sec2}, we recall some basic notions about exceptional groups and exceptional theta correspondences which are of our interest. Section \ref{sec3} is devoted to the local theory, while Section \ref{sec4} discusses the notion of relative cuspidality. In Section \ref{sec5}, we compute the necessary local L-factor, as a preparation for the global theory. In Section \ref{sec6}, we establish an Euler factorization of the global period integral and a global relative character identity. Finally, in Section \ref{sec7}, we construct the local relative characters of the degenerate Whittaker period associated with $F_4/\text{Spin}_9$ and carry out corresponding unramified computations.

\subsection*{Acknowledgment}

We would like to thank our advisor Wee Teck Gan for suggesting this problem, and for sending us an early draft of the preprint \cite{GG25}, as well as for his constant guidance and many useful comments. The authors also thank Edmund Karasiewicz, Spencer Leslie and Yi Shan for helpful discussions.
\section{Preliminaries}\label{sec2}
\subsection{The exceptional groups $F_4$ and $E_7$}
In this subsection, we recall some basic facts about the exceptional groups of interest.

Let $F$ be either a local field of characteristic $0$ or number field. Let $C$ be a composition algebra over $F$. It admits a conjugation $x\mapsto \bar{x}:C\rightarrow C$ and a quadratic norm form $N(x)=x\bar{x}$. We set a bilinear form associated to $N$ 
$$
B(x,y)=N(x+y)-N(x)-N(y),\text{ for }x,y\in C.
$$
and the trace $\text{Tr}(x)=x+\bar{x}$. When $\dim C = 8$, we have $C$ is isomorphic to the split octonian algebra $\mathbb{O}$ over $F$. Its automorphism group over $F$ is isomorphic to the split exceptional group of type $G_2$. We denote by 
$$
J=J(\mathbb{O})=\left\{ X= \begin{pmatrix}
a & z & \bar{y} \\
\bar{z} & b & x \\
y & \bar{x} & c
\end{pmatrix}\mid a,b,c\in F\text{ and }x,y,z\in \mathbb{O} \right\}
$$
the $27$-dimensional vector space consisting of $3\times 3$ Hermitian matrices whose entries are in $\mathbb{O}$. The set $J$ is equipped with the following multiplication
$$ X\circ Y = \frac{1}{2}(XY+YX), \text{ for }X,Y\in J,$$
where the multiplication on the right hand side is the usual matrix multiplication. We call $(J,+,\circ)$ the exceptional Jordan algebra. We set the trace of $X\in J$ to be $Tr(X)=a+b+c$, where $X=\begin{pmatrix}
a & z & \bar{y} \\
\bar{z} & b & x \\
y & \bar{x} & c
\end{pmatrix}$. This gives us the quadratic structure $q(X)=\text{Tr}(X\circ X)/2$, which is positive-definite. Its associated bilinear form is $(X,Y)=\text{Tr}(X\circ Y)$. We denote by $J_0$ the subspace of trace-$0$ elements in $J$. We also equip the exceptional Jordan algebra $J$ with the following cubic norm form given by the determinant map, i.e.
$$
N_J(X)= abc-aN(x)-bN(y)-cN(z)+\text{Tr}(xyz).
$$
We define the adjoint matrix of $X=\begin{pmatrix}
a & z & \bar{y} \\
\bar{z} & b & x \\
y & \bar{x} & c
\end{pmatrix}$ by
$$X^{\#}=\begin{pmatrix}
bc-N(x) & \overline{xy}-cz & zx-b\bar{y} \\
xy-c\bar{z} & ca-N(y) & \overline{yz}-ax \\
\overline{zx}-by & yz-a\bar{x} & ab-N(z)
\end{pmatrix}.$$
From this, we can define a symmetric cross product on $J$ by
$$
A \times B = ((A+B)^{\#}-A^{\#}-B^{\#})/2.
$$
The rank of $X$ is defined to be the degree of its minimal polynomial, i.e.
\begin{itemize}
    \item If $X=0$, then $X$ is of rank $0$.
    \item If $X\neq 0$ and $X^{\#}=0$, then $X$ is of rank $1$.
    \item If $X,\ X^{\#}\neq 0$ and $N_J(X)=0$, then $X$ is of rank $2$.
    \item Otherwise, $X$ is of rank $3$.
\end{itemize}

The group $\text{Aut}(J,\circ)$ consisting of automorphisms of $J$ over $F$ preserving the multiplication $\circ$ is a connected semisimple group over $F$ of type $F_4$. Its real point $F_4(\mathbb{R})$ is a subgroup of $\text{O}(J_\mathbb{R},q)$, which is to say it is compact. We denote this group by $F_4$. Moreover, we can identify $F_4$ with $\text{Aut}(J,N_J,I)$, the set of automorphisms of $J$ over $F$ preserving $N_J$ and the identity element $I$. As in \cite{Asc87}, maximal parabolic subgroups of $F_4$ can be described as stabilisers of singular subspaces of $J_0$. We set $P_1$ and $P_2$ to be maximal parabolic subgroups stabilising 1- and 2-dimensional singular subspaces of $J_0$, respectively. Let $M_J$ be the subgroup of $\text{Aut}(J)$ whose elements $m$ satisfy
$$
\det(m(A))=\lambda(m)\det(A),\,\forall A \in J,
$$
where $\lambda(m)\in F^\times$ is a similitude factor.

We consider the $56$-dimensional space $W_J=J\oplus F \oplus J \oplus F$ and equip it with the Freudenthal symplectic form and quartic form as follows.
For $w_i=(X_i,x_i,Y_i,y_i)\in W_J$, where $i\in \{1,2\}$, we set the following symplectic form
$$\langle w_1,w_2 \rangle = x_1y_2-x_2y_1+(X_1,Y_2)-(X_2,Y_1).$$
For $w=(X,x,Y,y)\in W_J$, we set the following quartic form
$$
Q(w)=(xy-(X,Y))^2+4xN_J(X)+4yN_J(Y)-4(X^{\#},Y^{\#}).
$$
We denote by $GE_7$ the group containing elements in $\text{Aut}(W_J)$ such that it preserves the above symplectic and quartic forms up to similitude $\nu:GE_7 \rightarrow \mathbb{G}_m$, i.e. $g\in \text{Aut}(W_J)$ satisfying 
$$
\langle gw_1,gw_2 \rangle = \nu(g) \langle w_1,w_2 \rangle \text{ and } Q(gw)=\nu(g)^2Q(w),
$$
for any $w,w_1,w_2 \in W_J$. The kernel of $\nu$ is a simply-connected Lie group of type $E_7$. It is of $F$-rank $3$.

The Siegel parabolic subgroup $P_J$ of $E_7$ is defined to be the subgroup stabilising the line generated by $(0,1,0,0)$. It has one Levi component whose elements also stabilise the line generated by $(0,0,0,1)$. This Levi subgroup is isomorphic to $M_J$ via the following action of $m\in M_J$ on $W_J$
$$
m(X,x,Y,y)=(m(X),\lambda(m)x,m^*(Y),\lambda(m)^{-1}y),
$$
where $m^*$ is the inverse adjoint of $m$. The unipotent radical $N_J$ of $P_J$ is abelian and isomorphic to $J$ as $F$-vector spaces. More concretely, this isomorphism is defined by the following action of $A$ on $W_J$:
$$
n(A)(X,x,Y,y)=(X+yA,x+(A,Y)+(A^{\#},X)+y\det(A),Y+2AX +yA^{\#},y).
$$
\subsection{Dual pairs $PGL_2\times F_4$ and $(SL_2\times SL_2)/\mu_2^\Delta \times \text{Spin}_9$ in $E_7$} We recall two dual pairs in $E_7$ for later use. 

We first consider the embedding $GL_2\times F_4 \hookrightarrow GE_7$. Let $J_0$ be the set of trace $0$ elements in $J$. This is an irreducible representation of $F_4$. Let $V_2=\langle u,v \rangle$ be the standard representation of $GL_2$. We define an action of $GL_2$ on $V_4 \cong F^4$ by
\begin{itemize}
    \item $\begin{pmatrix}
x &  \\
 & y
\end{pmatrix}(a,b,c,d)=(ya,x^2y^{-1}b,xc,x^{-1}y^2d)$
    \item $\begin{pmatrix}
1 & x \\
 & 1
\end{pmatrix}(a,b,c,d)=(a+xd,3x^2a+b+3xc+x^2d,2xa+c+x^2d,d)$
    \item $\begin{pmatrix}
0 & 1 \\
-1 & 0
\end{pmatrix}(a,b,c,d)=(-c,-d,a,b).$
\end{itemize}
We can identify $V_4\oplus V_2 \otimes J_0$ with $W_J$ via
$$
(a,b,c,d)+v\otimes X+u\otimes Y \mapsto(aI+X,b,cI+Y,d),
$$
whose inverse is
$$
(X,x,Y,y) \mapsto (\text{Tr}(X)/3,x,\text{Tr}(Y)/3,y)+v\otimes(X-\text{Tr}(X)I/3)+u\otimes(Y-\text{Tr}(Y)I/3).
$$
This identification gives us an embedding of $GL_2\times F_4$ into $GE_7$. By restricting to $E_7$ and taking modulo $\mu_2 = Z(E_7)$, we can deduce an embedding 
$$
PGL_2\times F_4 \hookrightarrow E_7.
$$
We set $e=\begin{pmatrix}
1 & 0 & 0 \\
0 & 0 & 0 \\
0 & 0 & 0
\end{pmatrix}.$ We can see that $F_4$ acts transitively on the set of rank $1$ elements in $J$ of trace $1$. The stabiliser of $e$ in $F_4$ is isomorphic to the group $\text{Spin}_9$ of type $B_4$. By \cite[Lemma 2.5.1]{Y25}, $\text{Spin}_9$ preserves
$$
\left\{\begin{pmatrix}
0 & 0 & 0 \\
0 & \xi & x \\
y & \bar{x} & -\xi
\end{pmatrix}\mid\xi\in F,\,x\in \mathbb{O}\right\}
\ \ \ \ \text{and}\ \ \ \ 
\left\{\begin{pmatrix}
0 & z & \bar{y} \\
\bar{z} & 0 & 0 \\
y & 0 & 0
\end{pmatrix}\mid y,z\in \mathbb{O}\right\}.
$$
It suffices to give an embedding of $(SL_2\times SL_2)/\mu_2^\Delta$ into $E_7$ whose image centralises $\text{Spin}_9$. Since the group $E_7$ we have constructed is of $F$-rank 3, we now construct a maximal $F$-split torus in this group. We define $3$ cocharacters of $E_7$ (as morphisms mapping $\mathbb{G}_m$ into $\text{Aut}(J)$):
\begin{itemize}
    \item $m_1(t):\begin{pmatrix}
a & z & \bar{y} \\
\bar{z} & b & x \\
y & \bar{x} & c
\end{pmatrix} \mapsto \begin{pmatrix}
t^2a & tz & t\bar{y} \\
t\bar{z} & b & x \\
ty & \bar{x} & c
\end{pmatrix}$,
    \item $m_2(t):\begin{pmatrix}
a & z & \bar{y} \\
\bar{z} & b & x \\
y & \bar{x} & c
\end{pmatrix} \mapsto \begin{pmatrix}
2a & tz & \bar{y} \\
t\bar{z} & t^2b & tx \\
y & t\bar{x} & c
\end{pmatrix}$,
    \item $m_3(t):\begin{pmatrix}
a & z & \bar{y} \\
\bar{z} & b & x \\
y & \bar{x} & c
\end{pmatrix} \mapsto \begin{pmatrix}
a & z & t\bar{y} \\
\bar{z} & b & tx \\
ty & t\bar{x} & tc
\end{pmatrix}$.
\end{itemize}
Then $m_i(t)$ are elements in $M_J$ satisfying $\lambda(m_i(t))=t^2$, for each $i=\overline{1,3}$. We have the following maximal split torus of $M_J$ and $E_7$
$$
T_1\times T_2 \times T_3/\mu_2^\Delta, 
$$
where $T_i=\{m_i(t)\mid\,t\in \mathbb{G}_m\}$ for $i=\overline{1,3}$ and $\mu_2^\Delta=\langle m_1(-1)m_2(-1)m_3(-1) \rangle$. The above split torus and unipotent elements
$$
n\left(\text{diag}(a,b,c)\right)\text{ and }\ n^\vee\left(\text{diag}(a,b,c)\right),
$$
for $a,b,c\in F$, generate a subgroup $\text{SL}_2\times \text{SL}_2 \times \text{SL}_2/\mu_2^\Delta$ in $E_7$. By \cite[Lemma 2.5.1]{Y25}, we can show that $\text{Spin}_9$ commutes with $m_1(t)$ and $m_2(t)m_3(t)$ for any $t\in \mathbb{G}_m$. Moreover, it also commutes with $n\left(\text{diag}(a,b,b)\right)$ and $n^\vee\left(\text{diag}(a,b,b)\right)$ for any $a,b\in F$. Therefore, we can construct a subgroup $\text{SL}_2\times \text{SL}_2/\mu_2^\Delta$ commuting with $\text{Spin}_9$ as desired.
\subsection{Measures} Let $F$ be a local field of characteristic $0$. We fix a nontrivial unitary additive character $\psi$ of $F$. We denote by $d_\psi x$ the self-dual additive Haar measure on $F$ with respect to $\psi$. Let $d^\times_\psi x=d_\psi x/|x|$ be a multiplicative Haar measure on $F^\times$. From now, since $\psi$ is fixed, we simply write the above measures as $dx$ and $d^\times x$.
\subsubsection{The group $PGL_2$} We consider the group $PGL_2$ over $F$. Let $B=TN$ be the upper triangular Borel subgroup, where
$$
T=\left\{\begin{pmatrix}
a & 0 \\
0 & 1
\end{pmatrix}\mid\ a\in F^\times\right\}\text{and }
N=\left\{\begin{pmatrix}
1 & b \\
0 & 1
\end{pmatrix}\mid\ b\in F\right\}.
$$
From the above identification, we can see that $T(F) \cong F^\times$ and $N(F)\cong F$. Therefore, we can define Haar measures on $T(F)$ and $N(F)$ via $dx$ and $d^\times x$. This gives us a right-invariant measure on $B(F)$. We also fix a Haar measure $dg$ on $PGL_2$, which determines a Plancherel measure $d\mu_{\text{PGL}_2}$ on the unitary dual $\widehat{\text{PGL}}_2$.

\subsubsection{The cone $\Omega$}\label{sec:ConeOmega} We set $\Omega$ to be the set of rank-1 matrices in the exceptional Jordan algebra and $\Omega_a\subseteq \Omega$ containing those with trace $a$, where $a\in F$. For any $a\neq 0$, we notice that $F_4$ acts transitively on $\Omega_a$ and the stabiliser of $e_a=\begin{pmatrix}
a & 0 & 0 \\
0 & 0 & 0 \\
0 & 0 & 0
\end{pmatrix}$ is isomorphic to $\text{Spin}_9$. Thus, we can identify $\Omega_a$ with $\text{Spin}_9\backslash F_4$. The Haar measures $d_\psi v$ on $\Omega$ and $d_\psi x$ on $F$ induce an $F_4$-invariant measure $|\omega_a|$ for each fiber $\Omega_a$, given $a\neq 0$. The measure $|\omega_a|$ is characterised by the following property: for any compactly supported smooth function $f$ on $\Omega \backslash \{0\}$, we have
$$
\int_\Omega f(v) d_\psi v =\int_{F^\times} \left(\int_{\Omega_a}f\cdot |\omega_a|\right)d_\psi a.
$$
\subsection{Exceptional Theta correspondences}
\subsubsection{The minimal representation of $E_7$}\label{sec:MinRepE7}
Since $P_J$ is a Siegel parabolic subgroup of $E_7$, we have an analog of the Schrodinger model for the minimal representation of $E_7$. As in \cite{DS99}, we realise the minimal representation $\Pi_\psi$, depending on a nonadditive character $\psi$ of $F$, on the space $L^2(\Omega)$ of square-integrable functions on $\Omega$. We have the following action of $\text{PGL}_2\times F_4$ on $\Pi_\psi$:
\begin{itemize}
    \item $(h\cdot f)(v)=f(h^{-1}\cdot v),$ for $h\in F_4$;
    \item $(t(a)\cdot f)(v)=|a|^6f(av)$, for $t(a)=\begin{pmatrix} a & 0 \\ 0 & 1 \end{pmatrix} \in T$;
    \item $(n(b)\cdot f)(v)=\psi(b\text{Tr}(v)) f(v)$, for $n(b)=\begin{pmatrix} 1 & b \\ 0 & 1 \end{pmatrix} \in N$;
    \item $(w\cdot f)(v)=\int_{F^\times}\int_\Omega f(x)\psi(r\cdot \langle v,x \rangle)|\eta(x)|dx \psi(r^{-1})|r|^3 d^\times r,$ for $w=\begin{pmatrix} 0 & 1 \\ 1 & 0 \end{pmatrix}$ is a Weyl group element of $\text{PGL}_2$.
\end{itemize}
We can consider the smooth representation $\Pi_\psi^\infty$ realised on the subspace $\mathcal{S}(\Omega)$ of Schwartz-Bruhat functions on $\Omega$. By a slight abuse of notation, we denote by $\Pi_\psi$ the minimal representation of $E_7$ in both smooth and $L^2$- settings.
\subsubsection{Smooth exceptional theta correspondence}
In this subsection, we recall the theory of exceptional theta correspondence for the dual pair $\text{PGL}_2\times F_4$ in the smooth setting.

For $\sigma \in \text{Irr}(\text{PGL}_2)$, the big theta lift of $\sigma$ to $F_4$ is
$$
\Theta_\psi(\sigma)=(\Pi_\psi \otimes \sigma^\vee)_{\text{PGL}_2},
$$
where the above subscript denotes the $\text{PGL}_2$-coinvariant space. By definition, we have the following $\text{PGL}_2$-invariant and $F_4$-equivariant projection map
$$
A_\sigma: \Pi_\psi \otimes \sigma^\vee \rightarrow \Theta_\psi(\sigma).
$$
This gives us a canonical $\text{PGL}_2\times F_4$-equivariant map
$$
\theta_\sigma:\Pi_\psi \rightarrow \sigma \boxtimes \Theta_\psi(\sigma).
$$
By Theorem 1.1 in \cite{KS23}, it follows that $\Theta_\psi(\sigma)$ is a nonzero finite-length representation of $F_4$ that has a unique irreducible quotient $\theta_\psi(\sigma)$. When $\sigma$ is tempered, $\Theta_\psi(\sigma)$ is irreducible and thus coincides with $\theta_\psi(\sigma)$. Moreover, if $\theta_\psi(\sigma_1)\cong \theta_\psi(\sigma_2)$, we have $\sigma_1 \cong \sigma_2$. The explicit description of $\Theta_\psi(\sigma)$ is given in \cite[Theorem 1.2]{KS23}.
\begin{thm}\label{2.1}
    Let $\sigma$ be an irreducible smooth representation of $\text{PGL}_2(F)$. Recall that $P_1$ and $P_2$ are maximal parabolic subgroups of $F_4$ defined in Section 2.1.
    \begin{enumerate}
        \item If $\sigma$ is a quotient of a principal series $\text{Ind}^{\text{PGL}_2}_{\bar{B}}(\chi)$, then $\Theta_\psi(\sigma)$ is a quotient of $\text{Ind}^{F_4}_{P_1}(\chi)$.
        \item If $\sigma$ is supercuspidal, then $\Theta_\psi(\sigma)=\theta_\psi(\sigma)$ is the unique irreducible quotient of $\text{Ind}^{F_4}_{P_2}(\chi\otimes|\det|^{3/2})$.
    \end{enumerate}
\end{thm}

Composing the maps $A_\sigma$ and $\theta_\sigma$ with the natural projection $\Theta_\psi(\sigma) \twoheadrightarrow \theta_\psi(\sigma)$, we have canonical equivariant maps
$$
A_\sigma: \Pi_\psi\otimes \sigma^\vee \rightarrow \theta_\psi(\sigma)
$$
and
$$
\theta_\sigma: \Pi_\psi \rightarrow \sigma \boxtimes \theta_\psi(\sigma),
$$
which are also denoted by the same symbols. Similarly, for $\pi \in \text{Irr}(F_4)$, we also have the following $F_4$-invariant and $\text{PGL}_2$-equivariant projection map 
$$
B_\pi: \Pi_\psi \otimes \pi^\vee \rightarrow \theta_\psi(\pi).
$$

We introduce the doubling zeta integral for later use. For $\Phi_1, \Phi_2 \in \mathcal{S}(\Omega)$ and $v_1,v_2\in \sigma$, the local doubling zeta integral is given by
$$
Z_\sigma(\Phi_1,\Phi_2,v_1,v_2)= \int_{\text{PGL}_2}\langle g\cdot \Phi_1,\Phi_2\rangle_\Omega \cdot \overline{\langle\sigma(g)\cdot v_1,v_2 \rangle_\sigma}dg,
$$
which converges when $\sigma$ is tempered. This defines a $(\text{PGL}_2\times \text{PGL}_2)$-invariant and $F_4$-invariant map 
$$
Z_\sigma: \Pi_\psi \otimes \bar{\Pi}_\psi \otimes \bar{\sigma} \otimes \sigma \rightarrow \mathbb{C}.
$$
Since $\bar{\sigma}\cong \sigma^\vee$, it follows that $Z_\sigma$ factors through the canonical projection map
$$
A_\sigma \otimes \overline{A_\sigma}: \Pi_\psi \otimes \bar{\Pi}_\psi \otimes \sigma^\vee \otimes \overline{\sigma^\vee} \rightarrow \Theta_\psi(\sigma)\otimes \overline{\Theta_\psi(\sigma)},
$$
which gives us
$$
Z_\sigma(\Phi_1,\Phi_2,v_1,v_2)=\langle A_\sigma(\Phi_1,v_1),A_\sigma(\Phi_2,v_2)\rangle_{\theta(\sigma)}
$$
for some Hermitian form $\langle-,-\rangle_{\theta(\sigma)}$ on $\Theta_\psi(\sigma)$. A reformulation of the above equality by $\theta_\sigma$ instead of $A_\sigma$ is
$$
\langle\theta_\sigma(\Phi_1),\theta_\sigma(\Phi_2)\rangle_{\sigma\boxtimes \theta(\sigma)}=\sum_{v\in \text{ONB}(\sigma)}Z_\sigma(\Phi_1,\Phi_2,v,v),
$$
where the above sum is over an orthonormal basis of $\sigma$.
\subsubsection{$L^2$-exceptional theta correspondence}\label{sec2.4.3} We consider the theta correspondence in the $L^2$-setting. As in \cite{GG14}, we have the following direct integral decomposition of $\text{PGL}_2\times F_4$-representations
$$
\Pi_\psi=L^2(\Omega)=\int_{\widehat{\text{PGL}_2}} \sigma\boxtimes \theta_\psi(\sigma) d\mu_{\text{PGL}_2}(\sigma),
$$
where $\mu_{\text{PGL}_2}$ is the Plancherel measure of $\text{PGL}_2(F)$ which depends on our fixed choice of measure on $\text{PGL}_2(F)$. By adapting \cite[Proposition 3.3.1]{S17} to the above decomposition, we give a spectral decomposition of $\langle-,-\rangle_\Omega$ as follows
$$
\langle\Phi_1,\Phi_2\rangle_\Omega = \int_{\widehat{\text{PGL}_2}}J^\theta_\sigma(\Phi_1,\Phi_2)d\mu_{\text{PGL}_2}(\sigma),
$$
where 
$$
J^\theta_\sigma(\Phi_1,\Phi_2)=\langle\theta_\sigma(\Phi_1),\theta_\sigma(\Phi_2)\rangle_{\sigma\boxtimes \theta(\sigma)}=\sum_{v\in \text{ONB}(\sigma)}Z_\sigma(\Phi_1,\Phi_2,v,v).
$$
Similar to \cite[Proposition 6.4]{GW18}, we have an explicit spectral decomposition of $L^2(\Omega_1,|\omega_1|)$ as follows
$$
L^2(\Omega_1,|\omega_1|)\simeq \int_{\widehat{\text{PGL}_2}}\dim \text{Hom}_N(\sigma,\psi)\cdot \theta_\psi(\sigma) d\mu_{\text{PGL}_2}(\sigma),
$$
which leads to the following commutative diagram
\[\begin{tikzcd}
                                                             & \mathcal{S}(\Omega) \arrow[ld, "rest"'] \arrow[rd, "\theta_\sigma"] &                                                                  \\
\mathcal{S}(\Omega_1) \arrow[rd, "\alpha_{\theta(\sigma)}"'] &                                                                     & \sigma\boxtimes\theta(\sigma) \arrow[ld, "{\ell_{\sigma,\psi}}"] \\
                                                             & \mathbb{C}_{\psi}\otimes \theta(\sigma)                             &                                                                 
\end{tikzcd}\]
where $\ell_{\sigma,\psi}$ is the $\psi$-Whittaker functional arising from the Whittaker-Plancherel theorem for $L^2(N,\psi\backslash \text{PGL}_2)$ and $\alpha_{\theta(\sigma)}$ is the morphism associated to the spectral decomposition of $L^2(\Omega_1,|\omega_1|)$.
\section{Relative character identity and transfer formula}\label{sec3}

In this section, we discuss the local theory. The aim is to give a conceptual, geometric definition of the transfer and relevant spaces of test functions, and establish the local relative character identities. We show that our definition of transfer via the exceptional theta correspondence essentially coincides with Sakellaridis' formula (Theorem \ref{thm3.5}). 

Retain the notation of \cite{GG25}. We also assume that the conductor of the additive character $\psi$ is $\O_F$. 

Recall the cone $\Omega$ from section \ref{sec:ConeOmega}, and recall $\text{Spin}_9\bs F_4 \cong \Omega_1\subseteq \Omega$ the subset of trace 1 elements in $\Omega$. Denote by $e\in \Omega_1$ the element $\text{diag}(1,0,0)$. 

\commentcustom{
\begin{lem} ($PGL_2$-action on $\Pi$)
Denote \[ t(a):=  \begin{pmatrix} a & 0 \\ 0 & 1 \end{pmatrix}, n(b):=  \begin{pmatrix}1 & b \\ 0 & 1 \end{pmatrix},   \] $w$ a Weyl group element of $PGL_2$ as in \cite{GG25}.
Then for $f\in \Pi$,
\begin{itemize}
    \item $(t(a)\cdot f)(v) = |a|^{6} f(av) $
    \item $(n(b) \cdot f)(v) = \psi(b Tr(v)) f(v)$
    \item $(w\cdot f)(v) = \int_{F^\times} \int_\Omega f(x) \psi(r\cdot \langle v,x \rangle ) \mathop{|\eta(x)|} dx \psi(r^{-1}) |r|^3 \mathop{d^\times r} $. 
\end{itemize}
\end{lem}
\begin{proof}
\cite[Theorem 1.1]{GG25}, \cite[Proposition 5.1]{MS97}
\end{proof}
}
    
\begin{lem}\label{lem:Submersion}
One has the submersion \begin{align*}
r : \Omega_1 &\rightarrow F \\ 
y &\mapsto \langle e,y\rangle
\end{align*}
away from the preimages of $0,1 \in F$. For $a\ne 0,1$, the fiber $r^{-1}(a)$ is a homogeneous space under $\text{Spin}_9$. So we may identify \[ \Omega_1^\heartsuit / \text{Spin}_9 \cong F^\heartsuit, \] where $F^\heartsuit := F \bs \{ 0, 1 \}$ and $\Omega_1^\heartsuit := r^{-1}(F^\heartsuit)$. 
\end{lem}
\begin{proof}
The exclusion of $a=0$ is clear. Now, noting that for $a\ne 0$, the fiber $r^{-1}(a)$ consists of those elements of $\Omega_1$ whose $(1,1)$-entry is $a$. For fixed $a\ne 0$, the $(2,1)$- and $(3,1)$-entries of elements of $\Omega_1$ form the 16-dimensional spin representation of $\text{Spin}_9$ under the action of $\text{Spin}_9$ on $\Omega_1$. In order that the trace be 1, one needs these two entries to have a corresponding norm, in the 16-dimensional spin representation, of $a(1-a)$ (cf. \cite[Section 2.2]{GG25}). Hence, upon excluding $a=1$, it remains to note that the vectors of fixed non-zero norm in the 16-dimensional spin representation of $\text{Spin}_9$, form a single orbit under the $\text{Spin}_9$ action. 

%(cf. also \url{https://math.ucr.edu/home/baez/octonions/node12.html#OP2})
\end{proof}

\begin{definition}\label{def:TransferMaps}
We define maps as follows:

\[\begin{tikzcd}
& \S_c(\Omega) \arrow[dl,"{p}"]  \arrow[rd, "{q=rest}"] &\\
\S(N,\psi \bs PGL_2)\arrow[d,"{\I}"] & & \S_c(\Omega_1) \arrow[d,"{\I^\prime}"]\\
\S(N,\psi \bs PGL_2)_{N,\psi}\arrow[d] & &\S_c(\Omega_1)_{\text{Spin}_9} \arrow[d] \\ 
C^\infty(F^\times)& &C^\infty(F^\heartsuit) \\ 
\end{tikzcd}\]

where here,

\begin{itemize}
    \item $p$ sends \[ p : f \mapsto \big( g\mapsto (g\cdot f)(e) \big) \]
    \item $q$ is the restriction map from $\Omega$ to $\Omega_1$;
    \item $\I$ is the $(N,\psi)$-orbital integral, defined by 
    \[ \I(p(f))(a) = \lim_{n\rightarrow\infty} \int_{U_n} (wt(a)n(b)\cdot f)(e) \overline{\psi(b)} \mathop{db} \]
    where $U_n:= \varpi^{-n} \O_F$. 
    The well-definedness of $\I$ will be shown in Lemma \ref{lem:OrbitalIntegralWellDefined} below. 
    \item $\I^\prime$ is the integration over fibers of the submersion $r$ (cf. Lemma \ref{lem:Submersion}).
    \item Finally, we identify $\S(N,\psi \bs PGL_2)_{N,\psi}$ and $\S_c(\Omega_1)_{\text{Spin}_9}$ with spaces of functions on $F^\times$ and $F^\heartsuit$ respectively, by regarding the former as $(N,\psi)$-invariant functions on the open Bruhat cell $NwB$ and hence on the torus $wT$, and the latter via the submersion of Lemma \ref{lem:Submersion}. 
\end{itemize}
\end{definition}

\begin{lem}\label{lem:OrbitalIntegralWellDefined}
The $(N,\psi)$-orbital integral $\I$, as given in Definition \ref{def:TransferMaps} above, is well-defined.  
\end{lem}
\begin{proof}
First, we recall the actions of $\text{PGL}_2\times F_4$ on the minimal representation $\Pi\cong L^2(\Omega)$, as written down in section \ref{sec:MinRepE7}. 

Then we have \begin{align*}
& \lim_{n\rightarrow\infty} \int_{U_n} (wt(a)n(b)\cdot f)(e) \overline{\psi(b)} \mathop{db} \\
&= \lim_{n\rightarrow\infty} \int_{U_n} \int_{F^\times} \int_\Omega (t(a)n(b)\cdot f)(v) \psi(r\cdot \langle v,e\rangle) \mathop{|\eta(v)|} \psi(r^{-1}) |r|^3 \mathop{d^\times r} \overline{\psi(b)} \mathop{db} \\
&= \lim_{n\rightarrow\infty} \int_{F^\times} \int_{U_n}  \int_\Omega (t(a)n(b)\cdot f)(v) \psi(r\cdot \langle v,e\rangle) \mathop{|\eta(v)|} \psi(r^{-1}) |r|^3  \overline{\psi(b)} \mathop{db} \mathop{d^\times r} \\
\end{align*}
where here we have interchanged the outer two integrals, by absolute convergence of the integral over $F^\times$ \cite[Proposition 5.2]{GG25},
\begin{align*}
&= |a|^{6} \lim_{n\rightarrow\infty} \int_{F^\times} \int_{U_n}  \int_{a^{-1} supp f} f(av) \psi(b Tr(av)) \psi(r\cdot \langle v,e\rangle) \mathop{|\eta(v)|}   \overline{\psi(b)} \mathop{db} \psi(r^{-1}) |r|^3 \mathop{d^\times r} \\
&= |a|^{-6} \lim_{n\rightarrow\infty} \int_{F^\times} \int_{U_n}  \int_{supp f} f(y) \psi(b Tr(y)) \psi(ra^{-1} \cdot \langle y,e\rangle) \mathop{|\eta(y)|}   \overline{\psi(b)} \mathop{db} \psi(r^{-1}) |r|^3 \mathop{d^\times r} \\
\end{align*}
where here we have made the substitution $y=av$, noting also the choice of measure on $\Omega$ (cf. \cite[Section 2.9]{GG25}),
\begin{align*}
&= |a|^{-6} \lim_{n\rightarrow\infty} \int_{F^\times}  \int_{supp f} f(y)  \psi(ra^{-1} \cdot \langle y,e\rangle) \big( \int_{U_n} \psi(b (Tr(y)-1))  \mathop{db}\big) \mathop{|\eta(y)|}    \psi(r^{-1}) |r|^3 \mathop{d^\times r} \\
&=  |a|^{-6} \lim_{n\rightarrow\infty}  \int_{F^\times} q^n \int_{supp f} f(y)  \psi(ra^{-1} \cdot \langle y,e\rangle) \big( \textbf{1}_{Tr(y) \in 1+\varpi^n \O_F} \big) \mathop{|\eta(y)|}    \psi(r^{-1}) |r|^3 \mathop{d^\times r} \\
&=  |a|^{-6}\int_{F^\times}  \int_{\Omega_1} f(y)  \psi(ra^{-1} \cdot \langle y,e\rangle) \mathop{|\eta(y)|}    \psi(r^{-1}) |r|^3 \mathop{d^\times r}, \\
\end{align*}
where here we have interchanged the limit and outer integral due to the boundedness result of \cite[Proposition 5.2(3)]{GG25}, applied to $q^n f(y)( \textbf{1}_{Tr(y) \in 1+\varpi^n \O_F})$ in place of $f$ and $a^{-1}e$ in place of $w$. More precisely, one needs the bound to be uniform in $n\rightarrow\infty$. But this bound is of the form $R_1(f)(w)$, which is an integration of $f$ over $\Omega^{sing}(0)$, cf. \cite[Remark 4.1]{GG25} and which follows by a similar argument as \cite[Proposition 3.16(2)]{GK23}. It follows that this bound is bounded as $n\rightarrow\infty$. 

Since this last expression is certainly well-defined, the lemma follows. 

\end{proof}

\commentcustom{
Although here we are integrating over $\Omega_1$, rather than the entire $\Omega$, one has a similar boundedness result for the integration over fixed $Tr(y)$.: see Lemma \ref{lem:BoundednessFixedTrace} below. 

\begin{lem}\label{lem:BoundednessFixedTrace}
Fix $f$, $a$ and $b$. Then
\[ \int_{\Omega_b} f(y)  \psi(ra^{-1} \cdot \langle y,e\rangle) \mathop{dy} |r|^4 \] 
is bounded for all $b$. Here $\Omega_b \subseteq \Omega$ is the set of (rank 1) trace $b$ elements. 
\end{lem}
\begin{proof}
Note that $supp f$ is compact. Now define $f^\prime(x) = f(bx/Tr(x)){\bf 1}_{supp f}$, and apply the bound in \cite[Proposition 5.2(3)]{GG25} to $f^\prime$. Observe that for fixed $b^\prime$, \[ \int_{\Omega_{b^\prime}} f^\prime(y)  \psi(ra^{-1} \cdot \langle y,e\rangle) \mathop{dy} |r|^4 \] is a bounded constant multiple of 
\[ \int_{\Omega_b} f(y)  \psi(ra^{-1} \cdot \langle y,e\rangle) \mathop{dy} |r|^4; \] the desired bound follows. 
\end{proof}
}

\begin{prop}\label{prop3.4}
The above maps  $p,q,\I,\I'$ give a well-defined transfer map \[ t_\psi : \S_c(\Omega_1)_{\text{Spin}_9} \rightarrow \S(N,\psi \bs PGL_2)_{N,\psi}.\]
This transfer map $t_\psi$ satisfies a relative character identity, similar to that stated in  \cite[Theorem 9.1]{GW18}: \[ \mathcal{B}_{\sigma,\ell_\sigma} \circ t_\psi = \mathcal{B}_{\theta_\psi(\sigma),\ell_{\theta_\psi(\sigma)}}. \]
Here, $\sigma$ is any tempered, $(N,\psi)$-generic representation of $PGL_2$ with (non-zero) theta lift $\theta_\psi(\sigma)$ to $\text{Spin}_9$; $\ell_\sigma $ and $\ell_{\theta_\psi(\sigma)}$ are the canonical functionals determined by the respective spectral decompositions of $N,\psi \bs PGL_2$ and $\text{Spin}_9\bs F_4$; and $\mathcal{B}_{\sigma,\ell_\sigma} , \mathcal{B}_{\theta_\psi(\sigma),\ell_{\theta_\psi(\sigma)}}$ are the relative characters defined similarly as in \cite[Section 7.4]{GW18}.
\end{prop}
\begin{proof}
The proof, which follows standard arguments, is completely analogous to that of \cite[Proposition 8.8]{GW18} and \cite[Theorem 9.1]{GW18}.
The only modifications are that one has $\S_c(\Omega)_{N,\psi} \cong \S_c(\Omega_1)$, which can be found in \cite[Section 8.1]{GG14}, and that one uses the spectral decomposition of $\text{Spin}_9\bs F_4$, which has already been shown in \cite[Theorem 16]{GG14}. 
\end{proof}

The following theorem, which gives an explicit formula for the transfer map $t_\psi$, is the main result of this section. 

\begin{thm}\label{thm3.5}
We regard $\S(N,\psi \bs PGL_2)_{N,\psi}$ and $\S_c(\Omega_1)_{\text{Spin}_9}$ as spaces of functions on $F^\times$ and $F^\heartsuit$ respectively (as in Definition \ref{def:TransferMaps}). Then an explicit formula for the transfer map $t_\psi$ is given by: \[ t_\psi(\phi)(a) =  |a|^{-3}\int_{F}  \int_{F} \phi(t)  \psi(tk^{-1}) \mathop{dt}    \psi(ka^{-1}) |k|^{-4} \mathop{dk}  \]
In other words, it is given by an iteration of two Fourier transforms.
\end{thm}
\begin{proof}

From the submersion in Lemma \ref{lem:Submersion}, we have
\begin{align*}
& |a|^{-6}\int_{F^\times}  \int_{\Omega_1} f(y)  \psi(ra^{-1} \cdot \langle y,e\rangle) \mathop{|\eta(y)|}    \psi(r^{-1}) |r|^3 \mathop{d^\times r} \\
 &= |a|^{-6}\int_{F}  \int_{F} \I^\prime(f)(t)  \psi(ra^{-1} \cdot t) \mathop{dt}    \psi(r^{-1}) |r|^2 dr 
\end{align*}
Now make the substitution $r=ak^{-1}$:
\begin{align*}
 = |a|^{-3}\int_{F}  \int_{F} \I^\prime(f)(t)  \psi(tk^{-1}) \mathop{dt}    \psi(ka^{-1}) |k|^{-4} \mathop{dk} 
\end{align*}
\end{proof}

\begin{rem}
We see that the above formula for $t_\psi$, which was obtained via the (exceptional) theta correspondence, essentially agrees with the formula given by Sakellaridis' results in \cite[Theorem 1.3]{S18}. 
\end{rem}

\section{Relative cuspidality}\label{sec4}
In this section, we discuss the notion of $X$-relative cuspidality. To be more precise, via the exceptional theta correspondence $\text{PGL}_2\times F_4\subset E_7$, we are able to construct $X$-relatively cuspidal representations of $F_4$ which are not supercuspidal in the usual sense.

We recall the definition of relative cuspidality in \cite{KT08}. Let $F$ be a nonarchimedean local field of residual characteristic not $2$. Let $G$ be a reductive group and $X=H\backslash G$ be a $G$-homogenenous space over $F$. Assume that we are in the case when $X(F)=G(F)/H(F)$. Let $(\pi,V)$ be an admissible representation of $G(F)$. For $\lambda \in \text{Hom}_H(\pi,\mathbb{C})$ and $v\in V$, let $\phi_{v,\lambda}$ denote the function on $X(F)$ defined by
$$
\phi_{v,\lambda}(x)=\lambda(\pi(x)v),\text{ for }x\in X(F).
$$
We call such functions $(X,\lambda)$-relative matrix coefficients of $\pi$.
\begin{definition}
    An $X$-distinguished representation $(\pi,V)$ is said to be $(X,\lambda)$-relatively cuspidal (for some $\lambda \in \text{Hom}_H(\pi,\mathbb{C})$) if all the $(X,\lambda)$-relative matrix coefficients of $\pi$ are compactly supported modulo $Z_G(F)H(F)$. A representation $\pi$ is said to be $X$-relatively cuspidal if $\pi$ is $(X,\lambda)$-relatively cuspidal for any $\lambda \in \text{Hom}_H(\pi,\mathbb{C})$.
\end{definition}
We consider the case when $X=\text{Spin}_9\backslash F_4$. Let $\sigma\in \text{Irr}_{\text{temp}}(\text{PGL}_2)$ be an irreducible tempered representation of $\text{PGL}_2(F)$. We denote $\pi=\theta_\psi(\sigma)$. Let $\ell_\sigma$ and $\ell_{\pi}$ be the canonical functional determined by the the spectral decompositions of $N,\psi\backslash \text{PGL}_2$ and $\text{Spin}_9\backslash F_4$. Let $v\in \sigma$ and $w\in \pi$. We denote by $\beta_\sigma(v)$ and $\beta_{\pi}(w)$ the functions on $N(F),\psi\backslash \text{PGL}_2(F)$ and $\text{Spin}_9(F)\backslash F_4(F)$ given by
$$
\beta_\sigma(v):g\mapsto \ell_\sigma(\sigma(g)v)\ \ \ \text{ and }\ \ \ 
\beta_{\pi}(w): g\mapsto \ell_\pi(\pi(g)w).
$$
We also recall the following projections
$$
A_\psi: \Pi_\psi \otimes \sigma^\vee \rightarrow\pi \ \ \ \text{and}\ \ \ 
B_\psi:\Pi \otimes \pi^\vee \rightarrow \sigma.
$$
We prove the following key proposition.
\begin{prop}\label{prop4.2}
    Given $\beta_\sigma$ and $\beta_\pi$ defined as above. Then $\beta_\sigma(v)$ is compactly supported on $(N(F),\psi)\backslash PGL_2(F)$ for any $v\in \sigma$ if and only if  $\beta_\pi(w)$ is compactly supported on $Spin_9(F)\backslash F_4(F)$ for any $w\in \pi$.
\end{prop}
\begin{proof}
    Assume $\beta_\sigma(v)$ is compactly supported for any $v\in \sigma$. It suffices to show that for any $f\in \Pi_\psi$ and $v\in \sigma$, we have
    $$
    g\mapsto \ell_\pi(\pi(g)A_\psi(f,v)): \text{Spin}_9(F)\backslash F_4(F) \rightarrow \mathbb{C}
    $$
    is compactly supported. By the weak Cartan decomposition for symmetric spaces (cf. \cite{BO07,DeSe11}), we have
    $$
    F_4(F)= \text{Spin}_9(F)\cdot A(F) \cdot K^\prime,
    $$
    for some 1-dimensional torus $A$ of $F_4$ and compact set $K^\prime$.
    We deduce to show that the function $a \mapsto \ell_\pi(\pi(a)A_\psi(f,v))$ is compactly supported on $A(F)$. Observe
    $$
    \ell_\pi(\pi(a)A_\psi(f,v))=\ell_\pi(A_\psi(\Pi_\psi(a)f,v))=\int_{N(F)\backslash \text{PGL}_2(F)} p(\Pi_\psi(a)f)(g)\cdot \overline{\beta_\sigma(v)(g)}\,dg
    $$
    $$
    =\int_{N(F)\backslash \text{PGL}_2(F)} (\Pi_\psi(g)f)(a^{-1}e)\cdot \overline{\beta_\sigma(v)(g)}\,dg
    =\int_{C} (\Pi_\psi(g)f)(a^{-1}e)\cdot \overline{\beta_\sigma(v)(g)}\,dg,
    $$
    for some compact set $C$ in $N(F)\backslash \text{PGL}_2(F)$ containing the support of $\beta_\sigma(v)$. Since $f$ and $\beta_\sigma(v)$ are smooth, it follows that the above integral is a linear combination of functions of the form
    $$
    a\in A(F) \mapsto f_i(a^{-1}e),\text{ for some }f_i\in C^\infty_c(\Omega).
    $$
    When $|a|\rightarrow 0$ or $|a|\rightarrow \infty$, it is clear that $a^{-1}e$ will escape the support of such $f_i$. In this case, we have $\ell_\pi(\pi(a)A_\psi(f,v))=\ell_\pi(A(\Pi_\psi(a)f,v))$ vanishes, which is to say the relative matrix coefficient $\ell_\pi(\pi(\cdot)A_\psi(f,v))$ is compactly supported.

    Conversely, suppose that $\beta_\pi(w)$ is compactly supported for any $w\in \pi$. 
    We want to show that for any $f\in \Pi_\psi$ and $w\in \sigma$, the function
    $$
    g\in \text{PGL}_2(F) \mapsto \ell_\sigma(\sigma(g)B_\psi(f,w))
    $$
    is compactly supported modulo $N(F)$. By the Iwasawa decomposition $\text{PGL}_2(F)=N(F)\cdot T(F)\cdot K$, it suffices to show that the support of the above function on $T(F)$ is compact. By the local relative character identity, we have
    $$
    \ell_\sigma(\sigma(t)B(f,w))=\ell_\sigma(B(\Pi_\psi(t)f,w))=\langle q(\Pi_\psi(t)f),\beta_\pi(w)\rangle_{X}
    $$
    $$
    =\int_{\text{Spin}_9(F)\backslash F_4(F)}(\Pi_\psi(t)f)(g^{-1}e)\cdot \overline{\beta_\pi(w)(g)}\,dg
    $$
    $$
    =|t|^{6}\int_{\text{Spin}_9(F)\backslash F_4(F)}f(tg^{-1}e)\cdot \overline{\beta_\pi(w)(g)}\,dg,
    $$
    where the last line follows from the description of the minimal representation $\Pi_\psi$ introduced in section \ref{sec:MinRepE7}. When $|t|\rightarrow \infty$, as $\beta_\pi(w)$ and $f$ are compactly supported, it follows that the above integral vanishes. We need to treat the case when $|t|$ sufficiently small separately. As in \cite[Section 5]{GG25}, for $f_0\in \mathcal{S}_c(\Omega)$, we recall the normalized Randon transform $\hat{\mathcal{R}}(f_0)$
    $$
    \hat{\mathcal{R}}(f_0)(w)=\int_\Omega f(v)\psi(\langle v,w\rangle_\Omega)|\eta(v)|.
    $$
    We consider the operator $\Phi=\Phi_1+\Phi_2$, where  
    $$
    \Phi_1(f_0)(w)=\int_F \hat{\mathcal{R}}(f_0)(rw)|r|^3d^\times r \ \text{ and }\ \Phi_2(f_0)(w)=\int_F \hat{\mathcal{R}}(f_0)(rw)(\psi(r^{-1})-1)|r|^3d^\times r.
    $$
    By \cite[Section 5]{GG25}, there exists $f_0\in \mathcal{S}_c(\Omega)$ such that $f=\Phi(f_0)$. This gives us 
    $$
    \int_{\text{Spin}_9(F)\backslash F_4(F)}f(tg^{-1}e)\cdot \overline{\beta_\pi(w)(g)}\,dg=\int_{\text{Spin}_9(F)\backslash F_4(F)}\Phi(f_0)(tg^{-1}e)\cdot \overline{\beta_\pi(w)(g)}\,dg
    $$
    $$
    =\int_{\text{Spin}_9(F)\backslash F_4(F)}\Phi_1(f_0)(tg^{-1}e)\cdot \overline{\beta_\pi(w)(g)}\,dg+\int_{\text{Spin}_9(F)\backslash F_4(F)}\Phi_2(f_0)(tg^{-1}e)\cdot \overline{\beta_\pi(w)(g)}\,dg.
    $$
When $|t|$ is sufficiently small, since $\beta_\pi(w)$ is compactly supported and $\Phi_2(f_0)$ is well-defined at the cone point, the second summand above is equal to
$$
\Phi_2(f_0)(0)\int_{\text{Spin}_9(F)\backslash F_4(F)}\overline{\beta_\pi(w)(g)}\,dg.
$$
The above integral is an $F_4$-invariant linear form of $\pi$. Since $\pi$ is not the trivial representation (otherwise $\beta_\pi(w)$ is not compactly supported), it follows that 
$$
\int_{\text{Spin}_9(F)\backslash F_4(F)}\overline{\beta_\pi(w)(g)}\,dg=0,\text{ for any }w\in \pi.
$$
We now consider the first summand
$$
\int_{\text{Spin}_9(F)\backslash F_4(F)}\Phi_1(f_0)(tg^{-1}e)\cdot \overline{\beta_\pi(w)(g)}\,dg.
$$
Since $\Phi_1(f_0)$ is homogeneous near $0$ (cf. \cite[Section 5]{GG25}), it follows that when $|t|$ is sufficiently small, the above integral is nonzero if and only if
\begin{equation}\label{eqn4.1}
    \int_{\text{Spin}_9(F)\backslash F_4(F)}\Phi_1(f_0)(g^{-1}e)\cdot \overline{\beta_\pi(w)(g)}\,dg \neq 0.
\end{equation}
By the discussion in page 20 in \cite{GG25}, the integral in (\ref{eqn4.1}) defines an $F_4$-invariant linear form of $\Pi_1\otimes\pi^\vee$, where $\Pi_1$ is the minimal representation of $E_6$. If the integral in (\ref{eqn4.1}) does not vanish, it follows that $\pi$ is the theta lift $\Theta^\pm$ of $\tau^\pm$ via the dual pair $\mu_2\times F_4 \subset E_6$, where $\tau^\pm$ is the trivial and nontrivial characters of $\mu_2$. In particular, as in \cite[Theorem 9.3]{KS23}, we have $\pi$ is an irreducible summand of the unique semisimple quotient of $i^{F_4}_Q(|\cdot|^{\frac{5}{2}}\circ\omega_4)$ (see \cite{KS23}). This gives us a contradiction since $\pi=\theta_\psi(\sigma)$ and $\sigma$ is tempered. Therefore 
$$
\int_{\text{Spin}_9(F)\backslash F_4(F)}\Phi_1(f_0)(g^{-1}e)\cdot \overline{\beta_\pi(w)(g)}\,dg = 0
$$
when $|t|$ is sufficiently small, which is to say $\ell_\sigma(\sigma(\cdot)B(f,w))$ is compactly supported. We have finished our proof.
\end{proof}
By Proposition \ref{prop4.2}, it follows that $\theta_\psi(\sigma)$ is $X$-relatively cuspidal if and only if $\sigma$ is $(N,\psi)\backslash \text{PGL}_2$-relatively cuspidal, and the latter is equivalent to $\sigma$ being a supercuspidal representations of $\text{PGL}_2$. We thus obtain the following theorem, characterising $X$-relatively cuspidal representations via the exceptional theta correspondence:
\begin{thm}\label{thm4.3}
    Let $\pi\in \text{Irr}(F_4)$ which is not $\Theta^\pm$. Then $\pi$ is $X$-relatively cuspidal if and only if it is the theta lift of a supercuspidal representation $\sigma$ of $\text{PGL}_2$.
\end{thm}

\section{Local L-factor}\label{sec5}
In this section, we study the local L-factor $L_X(s,-)$ associated to a certain spherical variety $X$. In general, this local L-factor has been computed in \cites{S08,S13}. By following the method in \cite{GW18}, when $X=\Omega_1=\text{Spin}_9\backslash F_4$, we compute $L_X(s,-)$ in terms of the local L-factor for the Whittaker variety $(N,\psi) \backslash PGL_2$ via an exceptional theta correspondence.

Let $F$ be a nonarchimedean local field of residual characteristic not 2. Let $\psi$ be an additive character whose conductor is $\mathcal{O}_F$ so that the associated measure $d_{\psi}x$ of $F$ gives $\mathcal{O}_F$ volume 1. We now define an integral structure on $X=\Omega_1$. We recall the Siegel parabolic subgroup $P_J=M_JN_J$. Let $\bar{N}_J$ be the opposite unipotent radical of $N_J$. Since our group $E_7$ is defined over $\mathbb{Z}$, this gives us an integral structure on $\bar{N}_J$. We fix a lattice $\Lambda\subset \bar{N}_J(F)$ such that $e=\text{diag}(1,0,0)$ lies in $\Lambda$. Let $K^\prime$ be the stabiliser of $\Lambda$ in $F_4(F)$. We set
$$
\Omega(n)=\Omega \cap (\varpi^n\Lambda \backslash \varpi^{n+1}\Lambda)
$$
so that 
$$
\Omega=\bigcup_{n\in \mathbb{Z}} \Omega(n).
$$
Then $\Omega(0)$ has volume $1$ with respect to the measure $d_\psi v$. Assume that the lattice $\Lambda$ endows $X=\Omega_1$ with the structure of a smooth scheme over $\mathcal{O}_F$. Then $e\in X(\mathcal{O_F})=X(F)\cap \Lambda$ and we have an $F_4$-equivariant isomorphism 
$$X\rightarrow \text{Spin}_9\backslash F_4$$
of smooth schemes over $\mathcal{O}_F$, which induces an isomorphism of their $\mathcal{O}_F$-points
$$X(\mathcal{O}_F) \cong \text{Spin}_9(\mathcal{O}_F)\backslash F_4(\mathcal{O}_F).$$
Moreover, the Haar measure $|\omega|$ on $X$ is endowed with $\frac{|\omega_{F_4}|}{|\omega_{\text{Spin}_9}|}$, where $\omega_{F_4}$ and $\omega_{\text{Spin}_9}$ are invariant differentials of top degree with nonzero reduction on the special fibers of smooth group schemes $F_4$ and $\text{Spin}_9$ over $\mathcal{O}_F$.
One has
$$\text{Vol}(X(\mathcal{O}_F);|\omega|)=\int_{X(\mathcal{O}_F)}|\omega|=q^{-\dim X} \cdot \frac{|F_4(\kappa_F)|}{|\text{Spin}_9(\kappa_F)|}=q^{-16} \cdot \frac{|F_4(\kappa_F)|}{|\text{Spin}_9(\kappa_F)|},$$
where $\kappa_F$ is the residue field of $F$ and $q$ is its cardinality.

We now study the local L-factor $L_X(s,-)$. Assume that $\sigma \in \widehat{\text{PGL}_2}$ is $K$-unramified, where $K=\text{PGL}_2(\mathcal{O}_F)$. We fix a $K$-spherical vector $v_0\in \sigma$ such that $\left\langle v_0,v_0 \right\rangle_\sigma=1$. Then $\theta_{\psi}(\sigma)$ is $K^\prime$-unramified and we fix a $K^\prime$-invariant vector $w_0\in \theta_{\psi}(\sigma)$ satisfying $\left\langle w_0,w_0 \right\rangle_{\theta_{\psi}(\sigma)}=1$. We define 
$$L^\#_X(\sigma)=|l_{\theta_{\psi}(\sigma)}(w_0)|^2.$$
We would like to determine $L^\#_X(\sigma)$ in terms of $|l_{\sigma}(v_0)|^2$. Let $\Phi_0$ be a spherical vector for the hyperspecial maximal compact subgroup $K^{\prime\prime}$ generating $\Pi_\psi$ (as a $P_J$-module). It is easy to see that $\Phi_0$ is $K\times K^\prime$-invariant. As discussed in \cite[Section 2.2]{Gan11}, it follows that $\Phi_0$ vanishes on $\Omega(n)$ if $n<0$ and $\Phi_0$ is constant on each $\Omega(n)$. Thus, $\Phi_0$ is determined by its values at $\varpi^n\chi_0$ for $n\geq 0$, where $\chi_0$ is a point in $\Omega(0)$. We normalised $\Phi_0$ so that it takes value $1$ on $\Omega(0)$. The main result in \cite{SW07} gives a formula for $\Phi_0$.
\begin{thm}\label{thm5.1}
    For any $n\geq 0$, we have
    $$
    \Phi_0(\varpi^n\chi_0)=\sum^n_{k=0}q^{3k}=\frac{q^{3(n+1)}-1}{q^3-1}.
    $$
\end{thm}
Under the map 
$$\theta_\sigma:\Pi_\psi \rightarrow \sigma \otimes \theta_\psi(\sigma),$$
we have
$$\theta_\sigma(\Phi_0)=c_\sigma \cdot v_0 \otimes w_0 \in \sigma \otimes \theta_\psi(\sigma),$$
for some nonzero constant $c_\sigma$.

Moreover, we have the basic functions
$$p(\Phi_0)=f_0\in \mathcal{S}(N,\psi\backslash \text{PGL}_2)\ \ \ \text{and}\ \ \ q(\Phi_0)=\phi_0\in S(X).$$
Observe $\phi_0=1_{X(\mathcal{O}_F)}$. We have 
$$\alpha_{\theta(\sigma)}(\phi_0)= \lambda_{\theta(\sigma)}\cdot w_0.$$
By taking inner product of both sides with $w_0$, we obtain
$$\lambda_{\theta(\sigma)}=\langle\alpha_{\theta(\sigma)}(\phi_0),w_0\rangle_{\theta(\sigma)}.$$
This gives us 
$$\langle\alpha_{\theta(\sigma)}(\phi_0),w_0\rangle_{\theta(\sigma)}=\langle\phi_0,\beta_{\theta(\sigma)}(w_0)\rangle_X=\beta_{\theta(\sigma)}(w_0)(v_1)\cdot \text{Vol}(X(\mathcal{O}_F);|\omega|).$$
Therefore $$\lambda_{\theta(\sigma)}=l_{\theta_{\psi}(\sigma)}(w_0)\cdot \text{Vol}(X(\mathcal{O}_F);|\omega|).$$
On the other hand, by the commutative diagram in section \ref{sec2.4.3}, it follows that
$$c_\sigma \cdot l_\sigma(v_0)\cdot w_0=l_\sigma(\theta_\sigma(\Phi_0))=\lambda_{\theta(\sigma)}\cdot w_0,$$
which gives us $$\lambda_{\theta(\sigma)}=c_\sigma \cdot l_\sigma(v_0).$$
Hence, we have 
$$|l_{\theta(\sigma)}(w_0)|=|l_{\sigma}(v_0)|\cdot |c_\sigma|\cdot \text{Vol}(X(\mathcal{O}_F);|\omega|)^{-1}.$$
It remains to determine $|c_\sigma|$. Observe
$$A_\sigma(\Phi_0,v_0)=\langle\theta_\sigma(\Phi_0),v_0\rangle_\sigma=c_\theta \cdot w_0.$$
By computing the inner product of both sides, it follows that
$$|c_\sigma|^2=\langle A_\sigma(\Phi_0,v_0),A_\sigma(\Phi_0,v_0)\rangle_{\theta(\sigma)}=Z_\sigma(\Phi_0,\Phi_0,v_0,v_0).$$
The following proposition gives us the neccessary numerical result.
\begin{prop}\label{prop:unramified} Given $\Phi_0$ and $v_0$ defined as above, we have  
$$|l_{\sigma}(v_0)|^2=\frac{\zeta_F(2)}{L(1,\sigma,Ad)}$$ and
$$Z_\sigma(\Phi_0,\Phi_0,v_0,v_0)=\text{Vol}(K;dg)\cdot \frac{L(\frac{11}{2},\sigma,std)\cdot L(\frac{5}{2},\sigma,std)}{\zeta_F(4) \cdot \zeta_F(8)},$$
where Ad and std are the adjoint and standard L-factor for $\text{PGL}_2$. Consequently, we have
$$L^\#_X(\sigma)=|l_{\theta_{\psi}(\sigma)}(w_0)|^2=\frac{L(\frac{11}{2},\sigma,std)\cdot L(\frac{5}{2},\sigma,std)}{\zeta_F(4) \cdot \zeta_F(8) \cdot L(1,\sigma,Ad)}$$
\end{prop}

\begin{proof}
    The calculation of $|l_{\sigma}(v_0)|$ was carried out in Proposition 2.14 in \cite{LM15}. As mentioned in \cite[page 42]{GW18}, we take the Haar measure associated to an invariant differential of top degree on $\text{PGL}_2$ over $\mathbb{Z}$. We have
    $$
    \text{Vol}(K;dg)=q^{-3}|\text{PGL}_2(\kappa_F)|=\zeta_F(2)^{-1}.
    $$
    It suffices compute the local doubling zeta integral. We denote by $\begin{pmatrix}
    \alpha_q & 0 \\
    0 & \alpha_q^{-1}
    \end{pmatrix}$ the Satake parameter of $\sigma$. We have the Cartan decomposition
    $$
    \text{PGL}_2(F)=\bigsqcup_{m\in \mathbb{Z}_{\geq 0}}K \begin{pmatrix}
    \varpi^m & 0 \\
    0 & 1
    \end{pmatrix}K,
    $$
    where $\varphi$ is the uniformizer of $\mathcal{O}_F$. Using the Cartan integral formula, we obtain
    $$
    Z_\sigma(\Phi_0,\Phi_0,v_0,v_0)=\int_{\text{PGL}_2(F)}\langle \Pi(g)\Phi_0,\Phi_0 \rangle_\Omega \cdot \overline{\langle \sigma(g)v_0,v_0 \rangle} dg
    $$
    $$
    =\sum_{m\geq 0} (q^{-1}+1)q^m \left\langle \Pi\left(\begin{pmatrix}
    \varpi^m & 0 \\
    0 & 1
    \end{pmatrix}\right)\Phi_0,\Phi_0 \right\rangle_\Omega \cdot \overline{\left\langle \sigma\left(\begin{pmatrix}
    \varpi^m & 0 \\
    0 & 1
    \end{pmatrix}\right)v_0,v_0 \right\rangle}.
    $$
    By the Macdonald formula (cf. \cite[Theorem 4.6.6]{B97}), it follows that
    $$
    \left\langle \sigma\left(\begin{pmatrix}
    \varpi^m & 0 \\
    0 & 1
    \end{pmatrix}\right)v_0,v_0 \right\rangle=\frac{1}{1+q^{-1}}q^{-m/2}\left[\alpha_q^m\frac{1-q^{-1}\alpha_q^{-2}}{1-\alpha_q^{-2}}+\alpha_q^{-m}\frac{1-q^{-1}\alpha_q^{2}}{1-\alpha_q^2}\right]
    $$
    $$
    =\frac{q^{-m/2}(q^{-1}\alpha_q^m+\alpha_q^{-m}-\alpha_q^{m+2}-q^{-1}\alpha_q^{-m+2})}{(1+q^{-1})(1-\alpha_q^2)}.
    $$
    By section \ref{sec:MinRepE7} and Theorem \ref{thm5.1}, we have
    $$
    \left\langle \Pi\left(\begin{pmatrix}
    \varpi^m & 0 \\
    0 & 1
    \end{pmatrix}\right)\Phi_0,\Phi_0 \right\rangle_\Omega=\int_{\Omega}|\varpi^m|^6 \Phi_0(\varpi^mv) \Phi_0(v)d_\psi v
    $$
    $$
    = q^{-6m}\sum_{n\geq 0} \frac{(q^{3(m+n+1)}-1)(q^{3(n+1)}-1)}{q^{17n}(q^3-1)^2}
    $$
    $$
    =\frac{q^{-3m+6}}{(q^3-1)^2(1-q^{-11})}-\frac{q^{-3m+3}}{(q^3-1)^2(1-q^{-14})}-\frac{q^{-6m+3}}{(q^3-1)^2(1-q^{-14})}+\frac{1}{(q^3-1)^2(1-q^{-17})}.
    $$
    Taking the above two formulas in $Z_\sigma(\Phi_0,\Phi_0,v_0,v_0)$ and by direct (but tedious) computation, we have
    $$
    Z_\sigma(\Phi_0,\Phi_0,v_0,v_0)=\frac{(1-q^{-4})(1-q^{-8})}{(1-q^{-5/2}\alpha_q)(1-q^{-5/2}\alpha_q^{-1})(1-q^{-11/2}\alpha_q)(1-q^{-11/2}\alpha_q^{-1})}
    $$
    $$
    =\frac{L(\frac{11}{2},\sigma,std)\cdot L(\frac{5}{2},\sigma,std)}{\zeta_F(4) \cdot \zeta_F(8)}.
    $$
\end{proof}

\section{Factorization of global periods}\label{sec6}
In this section, we study the factorization of global periods in terms of local ones. 

\subsection{Automorphic forms and global periods}

We define global analog terminologies of previous sections. Let $k$ be a number field with ring of adele $\mathbb{A}$. We fix a nontrivial unitary character $\psi:k\backslash \mathbb{A}\rightarrow S^1$. This character has a factorization $\psi = \underset{v}{\prod}\psi_v$, where $\psi_v$ is a nontrivial character of the local field $k_v$ for any place $v$ of $k$. Then for each place $v$, as in previous sections, one can define a self-dual Haar measure $d_{\psi_v}x$ on $k_v$. The product measure $dx=\underset{v}{\prod}d_{\psi_v}x$ gives us a measure on $\mathbb{A}$, which is a Tamagawa measure and does not depend on $\psi$.

Let $G$ be a reductive group over $k$ and $X=H\backslash G$ be a homogeneous $G$-variety over $k$. As mentioned in Section 12.1 in \cite{GW18}, by taking a product measure of local measures at every place, one can construct the Tamagawa measure $\omega_X$ of $X(\mathbb{A})$, which is equal to the quotient of the Tamagawa measures of $G(\mathbb{A})$ and $H(\mathbb{A})$.

Let $[G]$ be the quotient $G(k)\backslash G(\mathbb{A})$ and we equip it with the Tamagawa measure $dg$. We denote by $C^\infty_\text{mod}([G])$ the space of smooth functions of moderate growth on $[G]$. Let $\A(G)$ and $\A_\text{cusp}(G)$ be the space of automorphic forms and cusp forms on $G$. We have
$$\A_\text{cusp}(G) \subset \A(G) \subset C^\infty_\text{mod}([G])$$
as $G(\mathbb{A})$-modules.
The Petersson inner product $\langle -,- \rangle_G$ defines a pairing between $\A_\text{cusp}(G)$ and $C^\infty_\text{mod}([G])$. Hence, for any irreducible cuspidal representation $\Sigma \subset \A_\text{cusp}(G)$, we have the canonical projection map 
$$\text{pr}_\Sigma:C^\infty_\text{mod}([G]) \rightarrow \Sigma.$$
We consider the global $H$-period
$$P_H:\A_\text{cusp}(G) \rightarrow \mathbb{C}$$
defined by
$$P_H(\phi)=\int_{[H]}\phi(h)dh,$$
where $dh$ is the Tamagawa measure of $H(\mathbb{A})$. Let $P_{H,\Sigma}$ be the restriction of $P_H$ to $\Sigma$. We now define global analogs of $\alpha_\sigma$ and $\beta_\sigma$. Let
$$\begin{array}{ccccc}
\theta & : & C_{c}^{\infty}\left(X\left(\mathbb{A}\right)\right) & \longrightarrow & C_{\text{mod}}^{\infty}\left(\left[G\right]\right)\\
 &  & f & \mapsto & g\mapsto\underset{x\in X\left(k\right)}{\sum}f\left(x\cdot g\right).
\end{array}$$
Then one can define a composite map
$$\alpha_{\Sigma}^{\text{Aut}}:C_{c}^{\infty}\left(X\left(\mathbb{A}\right)\right)\overset{\theta}{\longrightarrow}C_{\text{mod}}^{\infty}\left(\left[G\right]\right)\overset{\text{pr}_{\Sigma}}{\longrightarrow}\Sigma.$$
We can see that 
$$\alpha_{\Sigma}^{\text{Aut}}(f)=\underset{\phi \in \text{ONB}(\Sigma)}{\sum} \langle \theta(f),\phi \rangle_G \cdot \phi.$$
On the other hand, we define
$$\begin{array}{ccccc}
\beta_{\Sigma}^{\text{Aut}} & : & \Sigma & \longrightarrow & C^{\infty}\left(X\left(\mathbb{A}\right)\right)\\
 &  & \phi & \mapsto & g\mapsto P_{H}\left(g\cdot\phi\right).
\end{array}$$
We recall Lemma 12.1 in \cite{GW18}.
\begin{lem}
    For $f\in C^\infty_c(X(\mathbb{A}))$ and $\phi \in \Sigma$, one has
    $$\langle \alpha_{\Sigma}^{\text{Aut}}(f),\phi \rangle_G = \langle f, \beta_{\Sigma}^{\text{Aut}}(\phi) \rangle_X.$$
\end{lem}
We now define a global analog of $J_\sigma$ and global relative characters. For $\phi_1,\phi_2 \in C^\infty_c(X(\mathbb{A}))$, we define
$$J^\text{Aut}_\Sigma (\phi_1,\phi_2)= \langle \alpha_{\Sigma}^{\text{Aut}}(\phi_1),\alpha_{\Sigma}^{\text{Aut}}(\phi_2) \rangle_\Sigma = \langle \beta_{\Sigma}^{\text{Aut}} \alpha_{\Sigma}^{\text{Aut}}(\phi_1), \phi_2 \rangle_X.$$
For $f\in C^\infty_c(X(\mathbb{A}))$, we define
$$B^\text{Aut}_\Sigma(f)=\beta_{\Sigma}^{\text{Aut}}(\alpha_{\Sigma}^{\text{Aut}}(f))(1)= \underset{\phi \in \text{ONB}(\Sigma)}{\sum}\langle f, \beta_{\Sigma}^{\text{Aut}}(\phi) \rangle_X \cdot P_H(\phi).$$

\subsection{Global exceptional theta correspondence for the dual pair $PGL_2 \times F_4$.}

We recall the global minimal representation for the dual pair $\text{PGL}_2\times F_4$ and related theta liftings. We make a remark that in this paper, we are working with the anisotropic form of $F_4$, namely $F_4(\mathbb{R})$ is compact. In this case, the archimedean theta correspondence for the dual pair $\text{PGL}_2\times F_4$ is completely determined in \cite{GS98}. For every place $v$, let $L_v=L\otimes_{\mathcal{O}_k}\mathcal{O}_{k_v}$. Then for almost all $v$, it is a self-dual lattice of volume $1$ with respect to $d_{\psi_v}$. Let $K^\prime_v$ be the stabilizer of $L_v$ in $F_4(k_v)$ and $\Phi_{0,v}=1_{L_v}\in S(\Omega_v)$, where $\Omega_v=\Omega \otimes_k k_v$. For each $v$, we have the minimal representation $\Pi_{\psi_v}$ of $\text{PGL}_2(k_v)\times F_4(k_v)$ realized on $\mathcal{S}(\Omega_v)$. Then we have the global minimal representation 
$$\Pi_\psi=\otimes^\prime_v\  \Pi_{\psi_v}$$
of $\text{PGL}_2(\mathbb{A})\times F_4(\mathbb{A})$ realized on
$$\mathcal{S}(\Omega_{\mathbb{A}})=\otimes^\prime_v\  \mathcal{S}(\Omega_v).$$
The minimal representation $\Pi_\psi$ has a canonical automorphic realization
$$\begin{array}{ccccc}
\theta & : & \Pi_{\psi} & \rightarrow & C_{\text{mod}}^{\infty}\left(\left[\text{PGL}_{2}\times F_{4}\right]\right)\\
 &  & \Phi & \mapsto & \left(g,h\right)\mapsto\underset{v\in\Omega_{k}}{\sum}\left(\Pi_{\psi}\left(g,h\right)\Phi\right)\left(v\right).
\end{array}$$
We are now able to define global theta lifting. Let $\Sigma \subset \mathcal{A}_\text{cusp}(\text{PGL}_2)$ be an irreducible cuspidal automorphic representation of $\text{PGL}_2(\mathbb{A})$. We define the global theta lift from $PGL_2$ to $F_4$ as follows. Let
$$\begin{array}{ccccc}
A^{\text{Aut}}_\Sigma & : & \Pi_{\psi}\otimes\bar{\Sigma} & \rightarrow & {\mathcal {A}}\left(F_{4}\right)\\
 &  & \Phi\otimes f & \mapsto & A^{\text{Aut}}_\Sigma(\Phi\otimes f)(h) = \int_{\left[\text{PGL}_{2}\right]} \theta\left(\Phi\right)\left(g,h\right)\cdot\overline{f\left(g\right)}dg.
\end{array}$$
We denote the image of the above map by $\Theta^{\text{Aut}}(\Sigma)$. If $\Theta^{\text{Aut}}(\Sigma)$ is cuspidal and nonzero, then by the Howe duality conjecture, one has $\Theta^{\text{Aut}}(\Sigma)$ is an irreducible cuspidal automorphic representation.

Conversely, let $\Sigma^\prime$ be an irreducible cuspidal automorphic representation of $F_4(\mathbb{A})$. We can define the global theta lift $B^\text{Aut}_{\Sigma^\prime}$ to $\text{PGL}_2$. By an interchange of order of integration, we have
$$\left\langle A_{\Sigma}^{\text{Aut}}\left(\Phi,f\right),\phi\right\rangle _{\Theta\left(\Sigma\right)}=\left\langle B_{\Theta\left(\Sigma\right)}^{\text{Aut}}\left(\Phi,\phi\right),f\right\rangle _{\Sigma},$$
for any $\Phi\in \Pi_\psi$ and $f\in \Sigma$ and $\phi \in \Theta^{Aut}(\Sigma)$.

\subsection{Global transfer of periods.}
Let $\Phi \in \Pi_\psi$ and $\phi \in \Sigma$, where $\Sigma$ is an irreducible cuspidal automorphic representation of $F_4(\mathbb{A})$. One can compute the $\psi$-Whittaker coefficient of $B^{\text{Aut}}_\Sigma (\Phi,\phi)$ in terms of $\text{Spin}_9$-period of $\phi$.
\begin{prop}
    For $\Phi \in \Pi_\psi$ and $\phi \in \Sigma$, one has
    $$P_{N,\psi} (B^{\text{Aut}}_\Sigma (\Phi,\phi)) = \int_{\text{Spin}_9(\mathbb{A})\backslash F_4(\mathbb{A})} \Phi(g^{-1}v_1) \cdot P_{\text{Spin}_9}(\phi)(g) \frac{dg}{dh} = \langle \Phi|_X, \beta^{Aut}_\Sigma(\phi) \rangle_{X_\mathbb{A}}.$$
    In particular, a cuspidal representation $\Sigma$ of $F_4$ has nonzero $\text{Spin}_9$-period if and only if its global theta lift to $\text{PGL}_2$ is globally $\psi$-generic.
\end{prop}
\begin{proof}
    We have $$P_{N,\psi}\left(B_{\Pi}^{Aut}\left(\Phi,\phi\right)\right)=\int_{\left[N\right]}\left(\int_{\left[F_{4}\right]}\theta\left(\Phi\right)\left(ng\right)\cdot\overline{\phi\left(g\right)}dg\right)\overline{\psi\left(n\right)}dn$$    $$=\int_{\left[F_{4}\right]}\left(\int_{\left[N\right]}\theta\left(\Phi\right)\left(ng\right)\cdot\overline{\psi\left(n\right)}dn\right)\overline{\phi\left(g\right)}dg,$$
    note that $[N]$ is compact.The inner integral is equal to 
    $$\underset{v\in\Omega\left(k\right)}{\sum}\int_{\left[N\right]}\left(\Pi_{\psi}\left(ng\right)\Phi\right)\left(v\right)\cdot\overline{\psi\left(n\right)}dn=\underset{v\in\Omega\left(k\right)}{\sum}\int_{\left[N\right]}\left(\Pi_{\psi}\left(g\right)\Phi\right)\left(v\right)\cdot\psi\left(n\text{Tr}\left(a\right)\right)\overline{\psi\left(n\right)}dn$$
    $$=\underset{v\in\Omega_{1}\left(k\right)}{\sum}\left(\Pi_{\psi}\left(g\right)\Phi\right)\left(v\right)=\theta\left(\Phi|_{X}\right)\left(g\right).$$
    On the other hand, observe 
    $$\langle\Phi|_{X},\beta_{\Sigma}^{Aut}(\phi)\rangle_{X_{\mathbb{A}}}=\langle\alpha_{\Sigma}^{Aut}(\Phi|_{X}),\phi\rangle_{F_{4}(\mathbb{A})}=\int_{\left[F_{4}\right]}\theta\left(\Phi|_{X}\right)\left(g\right)\cdot\overline{\phi\left(g\right)}dg,$$
    which is to say 
    $$P_{N,\psi}\left(B_{\Pi}^{Aut}\left(\Phi,\phi\right)\right)=\langle\Phi|_{X},\beta_{\Sigma}^{Aut}(\phi)\rangle_{X_{\mathbb{A}}}.$$
    
\end{proof}
\subsection{Factorization of global periods} In this subsection, we obtain a factorization of global periods. We first introduce some adelic periods, which are (regularized) Euler products of local periods.

Suppose $\Sigma = \otimes^\prime_v \Sigma_v$ be an irreducible tempered $\psi$-generic cuspidal automorphic representation of $\text{PGL}_2(\mathbb(A)$ and $\Sigma^\prime=\Theta^{\text{Aut}}(\Sigma)$ is a nonzero irreducible cuspidal representation of $F_4(\mathbb{A})$. Then $\Sigma^\prime$ is globally $\text{Spin}_9$-distinguished. On the other hand, by the main results in \cite{GS98} and \cite{KS23}, for each place $v$, we have the local big theta lift 
$\Theta_{\psi_v}(\Sigma_v) := (\Pi_{\psi_v} \otimes \Sigma_v^\vee)_{\text{PGL}_2(k_v)}$ and its unique irreducible quotient $\Sigma^\prime_v=\theta_{\psi_v}(\Sigma_v)$. We equip $\Sigma^\prime_v$ with the inner product arising from the local doubling zeta integral. Let $w_{0,v}\in \Sigma^\prime_v$ be a $K^\prime_v$-spherical unit vector for almost all $v$. We form the abstract theta lift
$$\Theta^{\mathbb{A}}(\Sigma)=\otimes^\prime_v\theta_{\psi_v}(\Sigma_v),$$
where the restricted tensor product is with respect to the family of unit vectors $w_{0,v}$. The abstract theta lift inherits a unitary structure from that of its local components. We fix an isomorphism
$$\Theta^{Aut}(\Sigma) \cong \Theta^{\mathbb{A}}(\Sigma)$$
so that 
$$\langle -,- \rangle_{\Theta^{Aut}(\Sigma)}=\underset{v}{\prod}\langle -,- \rangle_{\Sigma^\prime_v},$$
where the inner product on the LHS is defined by the Petersson inner product.

We define global analogs of various local objects by normalizing local factors. We first define a global analog of $\ell_{\Sigma_v}$. We set   
$$\ell_{\Sigma_{v}}^{\flat}=\lambda_{1,v}^{-1}\cdot\ell_{\Sigma_{v}},$$
with $\left|\lambda_{1,v}\right|^{2}=\frac{\zeta_{k_v}(2)}{L(1,\Sigma_v,Ad)}$
and $\ell_{\Sigma_{v}}^{\flat}(v_{0,v})=1$ for almost all $v$. Let 
$$\ell_\Sigma^{\mathbb{A}}=\left(\frac{|\zeta_k(2)|}{|L(1,\Sigma,Ad)|} \right)^{1/2} \cdot \underset{v}{\prod}\ell_{\Sigma_{v}}^{\flat},$$
where the global L-value is defined by meromorphic continuation of the global L-function. Then for any $u=\otimes_v u_v \in \Sigma$, with $u_v=v_{0,v}$ for almost all $v$, we have $\ell_\Sigma^{\mathbb{A}}(u)$ converges.

Next, we set
$$A_{\Sigma_{v}}^{\flat}=\lambda_{2,v}^{-1}\cdot A_{\Sigma_{v}},$$
with $\left|\lambda_{2,v}\right|^{2}=\frac{L(\frac{11}{2},\Sigma_v,std)\cdot L(\frac{5}{2},\Sigma_v,std)}{\zeta_{k_v}(4) \cdot \zeta_{k_v}(10)}$
and $A_{\Sigma_{v}}^{\flat}(\Phi_{0,v},v_{0,v})=w_{0,v}$, for almost all $v$. We define 
$$A_{\Sigma}^{\mathbb{A}}=\left(\frac{|L(\frac{11}{2},\Sigma,std)\cdot L(\frac{5}{2},\Sigma,std)|}{|\zeta_k(4) \cdot \zeta_k(8)|} \right)^{1/2} \cdot \underset{v}{\prod} A_{\Sigma_{v}}^{\flat}.$$

Let 
$$\ell_{\Sigma^\prime_{v}}^{\flat}=\lambda_{3,v}^{-1}\cdot \ell_{\Sigma_{v}},$$
with $\left|\lambda_{3,v}\right|^{2}= \frac{L(\frac{11}{2},\Sigma_v,std)\cdot L(\frac{5}{2},\Sigma_v,std)}{\zeta_F(4) \cdot \zeta_F(8) \cdot L(1,\Sigma_v,Ad)}$ and $\ell_{\Sigma^\prime_{v}}^{\flat}(w_{0,v})=1$ for almost all $v$. Then we set
$$\ell_{\Theta\left(\Sigma\right)}^{\mathbb{A}}=\left(\frac{\left|L(\frac{11}{2},\Sigma,std)\cdot L(\frac{5}{2},\Sigma,std)\right|}{\left|\zeta_k(4) \cdot \zeta_k(8) \cdot L\left(1,\Sigma,Ad\right)\right|} \right)^{1/2}\cdot\underset{v}{\prod}\ell_{\Theta\left(\Sigma_{v}\right)}^{\flat}.$$
We can now compare various adelic period maps. Since $\ell^{\mathbb{A}}_\Sigma$ and $P_{\sigma,N,\psi}$ are nonzero elements in the $1$-dimensional space $\text{Hom}_{N(\mathbb{A})}(\Sigma,\psi)$, there exists a nonzero complex number $c(\Sigma)$ such that
$$P_{\Sigma,N,\psi}=c(\Sigma) \cdot \ell^{\mathbb{A}}_\Sigma.$$
Then we have 
$$\beta_{\Sigma}^{\text{Aut}}=c(\Sigma) \cdot \beta_{\Sigma}^{\mathbb{A}}.$$
Similarly, there exists a nonzero complex number $c(\Theta(\Sigma))$ such that
$$P_{\Theta(\Sigma),\text{Spin}_9}=c(\Theta(\Sigma)) \cdot \ell^{\mathbb{A}}_{\Theta(\Sigma)}$$
and
$$\beta_{\Theta(\Sigma)}^{\text{Aut}}=c(\Theta(\Sigma)) \cdot \beta_{\Theta(\Sigma)}^{\mathbb{A}}.$$
Moreover, we have $a(\Sigma)$ and $b(\Sigma)$ in $\mathbb{C}^\times$ such that
$$A_\Sigma^{Aut}=a(\Sigma)\cdot A^{\mathbb{A}}_\Sigma \ \ \text{and}\ \ B_{\Theta(\Sigma)}^{Aut}=b(\Sigma)\cdot B^{\mathbb{A}}_{\Theta(\Sigma)}.$$
We determine the constant $|c(\Theta(\Sigma))|^2$ by relating it to other constants $a(\Sigma)$, $b(\Sigma)$ and $c(\Sigma)$. We state an analog of Proposition 12.3 in \cite{GW18}.
\begin{prop}
    We have
    $$\overline{c(\Theta(\Sigma))}=c(\Sigma) \cdot b(\Sigma).$$
\end{prop}
There remains to compute $c(\Sigma)$ and $b(\Sigma)$.
\begin{prop}
    We have
    $$a(\Sigma)=b(\Sigma)\ \ \ \text{and}\ \ \ c(\Sigma)=1.$$
\end{prop}
\begin{proof}
$a(\Sigma)=b(\Sigma)$ follows from Proposition 12.4 in \cite{GW18}. The computation for $c(\Sigma)$ follows from Section 6.1 in \cite{LM15}, noting that $\text{PGL}_2$ is isomorphic to $\text{SO}_3$.
\end{proof}
We now show that $|a(\Sigma)|=1$, which is to say $|c(\Theta(\Sigma))|=|c(\Sigma)|=1$.
\begin{prop}
    For $\Phi \in \Pi_\psi$ and $f\in \Sigma$, we have 
    $$P_{\text{Spin}_9}(A^{Aut}_\Sigma(\Phi,f))=\langle p(\Phi), \beta^{Aut}_\Sigma(f) \rangle_{N\backslash PGL_2}.$$
    Consequently, we have $|a(\Sigma)|=1$.
\end{prop}
\begin{proof}
    We consider the following see-saw diagram
    \[\begin{tikzcd}
	{(\text{SL}_{2}\times \text{SL}_{2})/\mu_2^\Delta} && {F_4} \\
	\\
	{\text{PGL}_2} && {\text{Spin}_9}
	\arrow[no head, from=1-1, to=3-3]
	\arrow[no head, from=3-3, to=1-3]
	\arrow[no head, from=3-1, to=1-3]
	\arrow[no head, from=1-1, to=3-1]
    \end{tikzcd}\]
    which gives us a global see-saw identity. To be more precise, we have 
    $$P_{\text{Spin}_9}(A^{Aut}_\Sigma(\Phi,f))=\int_{[\text{PGL}_2]}\overline{f(g)}\cdot I(\Phi)(g)dg,$$
    where $I(\Phi)$ is the theta integral defined by
    $$I(\Phi)(g):=\int_{[\text{Spin}_9]}\theta(\Phi)(gh)dh.$$
    Note that $I(\Phi)$ belongs to the global theta lift of the trivial representation of $\text{Spin}_9$ to $(\text{SL}_{2}\times \text{SL}_{2})/\mu_2^\Delta$. By Theorem 6.1.3 in \cite{Y25}, which is an exceptional Siegel-Weil formula for $\text{Spin}_9\times (\text{SL}_{2}\times \text{SL}_{2})/\mu_2^\Delta$, we have
    $$I(\Phi) = E(\text{Res}(\phi_\Phi)),$$
    where $\phi_\Phi(g)= (g\cdot \Phi)(0)$ as a function on $E_7(\mathbb{A})$, and $\text{Res}(\phi_\Phi)$ is its restriction to a function on the adelic points of $(\text{SL}_{2}\times \text{SL}_{2})/\mu_2^\Delta$. This gives us 
    $$P_{\text{Spin}_9}(A^{Aut}_\Sigma(\Phi,f))=\int_{[\text{PGL}_2]}\overline{f(g)}\cdot E(\text{Res}(\phi_\Phi))(g)dg.$$
    By unfolding the RHS as in Section 6.2 in \cite{Y25}, we have
    $$P_{\text{Spin}_9}(A^{Aut}_\Sigma(\Phi,f))=\int_{N(\mathbb{A})\backslash\text{PGL}_2(\mathbb{A})}\overline{P_{N,\psi}(g\cdot f)}\cdot p(\Phi)(g)dg.$$
    By combining the above result with its local counterpart, we obtain
    $$\overline{c(\Sigma)}=c(\Theta(\Sigma))\cdot a(\Sigma).$$
    Proposition 3.3 and Proposition 3.4 give us
    $$\overline{c(\Sigma)}=\overline{c(\Sigma)}\cdot \overline{b(\Sigma)}\cdot a(\Sigma)=\overline{c(\Sigma)}\cdot |a(\Sigma)|^2,$$
    which is to say 
    $$|a(\Sigma)|=1.$$
\end{proof}
The above discussion gives a proof for a factorization of global periods, which is stated in the following theorem.
\begin{thm}\label{thm6.6}
    Let $\Sigma$ be a globally $\psi$-generic cuspidal automorphic representation of $\text{PGL}_2$ such that $\Sigma^\prime = \Theta(\Sigma) \subset \mathcal{A}_\text{cusp}(F_4)$. Then $$\left|P_{\Sigma^{\prime},\text{Spin}_{9}}\left(\phi\right)\right|^{2}=\left|\ell_{\Sigma^{\prime}}^{\mathbb{A}}\left(\phi\right)\right|^{2},$$
    for any $\phi \in \Sigma^\prime$.
\end{thm}
\subsection{Global relative character identity}
In this section, we obtain the global analog of the relative character identity. We consider the following diagram
\[\begin{tikzcd}
& \S(\Omega_{\mathbb{A}}) \arrow[dl,"{p}"]  \arrow[rd, "{q}"] &\\
\S(N(\mathbb{A}),\psi \bs PGL_2(\mathbb{A})) & & \S(\Omega_{1,\mathbb{A}})  \\ 
\end{tikzcd}\]
The space $\mathcal{S}(N(\mathbb{A}),\psi \backslash \text{PGL}_2(\mathbb{A}))$ is the restricted tensor product of the local spaces $\mathcal{S}(N(k_v),\psi \backslash \text{PGL}_2(k_v))$ of test functions with respect to the family of basic functions ${f_{0,v}}$. Similarly, the space $\mathcal{S}(\Omega_{1,\mathbb{A}})$ is the restricted tensor product of $S(\Omega_{1,k_v})$ (which is $C^\infty_c(\Omega_{1,k_v})$ at finite places) with respect to the family of basic functions ${\phi_{0,v}}$. We say that $f\in \mathcal{S}(N(\mathbb{A}),\psi \backslash \text{PGL}_2(\mathbb{A}))$ and $\phi \in \mathcal{S}(\Omega_{1,\mathbb{A}})$ are transfers of each other if there exists $\Phi \in S(\Omega_{\mathbb{A}})$ such that $f=p(\Phi)$ and $\phi=q(\Phi)$. The fundamental lemma shows that every $f$ has a transfer $\phi$ and vice versa.

As discussed in Section 12.16 of \cite{GW18}, we need to ensure that the map $\alpha_\Sigma^\text{Aut}=\text{pr}_\Sigma \circ \theta$ can be defined on $\mathcal{S}(N(\mathbb{A}),\psi \backslash \text{PGL}_2(\mathbb{A}))$. Let $f= p(\Phi) = \otimes_v^\prime f_v\in \mathcal{S}(N(\mathbb{A}),\psi \backslash \text{PGL}_2(\mathbb{A}))$ for some $\Phi \in S(\Omega_\mathbb{A})$. We have 
$$\theta(f)(g)= \underset{\gamma \in N(k)\backslash \text{PGL}_2(k)}{\sum}f(\gamma g),\ \ \text{for }g\in \text{PGL}_2(\mathbb{A}).$$
We need to prove the following sum is convergent
$$\underset{\gamma \in B(k)\backslash \text{PGL}_2(k)}{\sum}\ \underset{\lambda \in k^\times}{\sum}|f(t(\lambda)\gamma g)| = \underset{\gamma \in B(k)\backslash \text{PGL}_2(k)}{\sum}\ \underset{\lambda \in k^\times}{\sum} |(\gamma g \cdot \Phi)(\lambda \cdot e_1)|.$$
We denote the inner sum by $F_\Phi(g)$, which is convergent since $\Phi$ is a Schwartz function on $\Omega_{\mathbb{A}}$. We need to show the following sum converges
$$\underset{\gamma \in B(k)\backslash \text{PGL}_2(k)}{\sum} F_\Phi(\gamma g)$$
This follows verbatim from the discussion in page 59 in \cite{GW18}, noting that 
$$F_\Phi (t(a)k) = |a|_{\mathbb{A}}^6 \cdot \underset{\lambda \in k^\times}{\sum}|(k\cdot \Phi)(a\lambda \cdot e_1)|,$$
for any $a\in \mathbb{A}^\times$ and $k \in K = \prod_v K_v$.
Hence, the map $\theta$ defines a $\text{PGL}_2(\mathbb{A})$-equivariant map
$$\theta:\mathcal{S}(N(\mathbb{A}),\psi \backslash \text{PGL}_2(\mathbb{A})) \longrightarrow C^\infty_\text{mod}([\text{PGL}_2]). $$
The map $\alpha^\text{Aut}_\Sigma= \text{pr}_\Sigma \circ \theta$ is well-defined and we still have the adjunction formula
$$\langle \alpha^\text{Aut}_\Sigma(f),\phi \rangle_\Sigma = \langle f,\beta^\text{Aut}_\Sigma (\phi) \rangle_{N(\mathbb{A})\backslash \text{PGL}_2(\mathbb{A})},$$
for any $f\in \mathcal{S}(N(\mathbb{A}),\psi \backslash \text{PGL}_2(\mathbb{A}))$ and $\phi \in \Sigma$.
We have the following global relative character identity.
\begin{thm}\label{thm6.7}
    If $f\in \mathcal{S}(N(\mathbb{A}),\psi \backslash \text{PGL}_2(\mathbb{A}))$ and $\phi \in \mathcal{S}(\Omega_{1,\mathbb{A}})$ are transfer of each other, then for a cuspidal representation $\Sigma$ of $\text{PGL}_2$ with cuspidal theta lift $\Theta^\text{Aut}(\Sigma)$ on $F_4$, one has $$\mathcal{B}^\text{Aut}_\Sigma(f)=\mathcal{B}^\text{Aut}_{\Theta(\Sigma)}(\phi).$$
\end{thm}
\begin{proof}
    We have 
    $$\beta_\Sigma^{Aut}=c(\Sigma)\cdot \beta_\Sigma^{\mathbb{A}}.$$
    Applying the adjunction formula for $(\alpha_\Sigma^{\text{Aut}},\beta_\Sigma^{\text{Aut}})$ and $(\alpha_\Sigma^{\mathbb{A}},\beta_\Sigma^{\mathbb{A}})$, we obtain
    $$\alpha_\Sigma^\text{Aut}=\overline{c(\Sigma)}\cdot \alpha_\Sigma^\mathbb{A}.$$
    Therefore
    $$\mathcal{B}^\text{Aut}_\Sigma(f)=(\beta_\Sigma^{\text{Aut}}\alpha_\Sigma^{\text{Aut}}(f))(1)=|c(\Sigma)|^2 \cdot (\beta_\Sigma^{\mathbb{A}}\alpha_\Sigma^{\mathbb{A}}(f))(1) = |c(\Sigma)|^2 \cdot \mathcal{B}^\mathbb{A}_\Sigma(f).$$
    Similarly, we have 
    $$\mathcal{B}^\text{Aut}_{\Theta(\Sigma)}(\phi)= |c(\Theta(\Sigma))|^2 \cdot \mathcal{B}^\mathbb{A}_\Sigma(\phi).$$
    Since $|c(\Sigma)|=|c(\Theta(\Sigma))|$, the global relative character identity is deduced from the local relative character identity (cf. Proposition \ref{prop3.4}).
\end{proof}

\section{Degenerate Whittaker period of $F_4/\text{Spin}_9$}\label{sec7}
Let $F$ be a nonarchimedean local field of characteristic $0$. In this section, we construct the local relative characters of the degenerate Whittaker period associated with $X=\text{Spin}_9\backslash F_4$ over $F$. This construction provides the refinement of \cite[Conjecture 1.6]{MWZ26} as mentioned in \cite[Remark 1.7]{MWZ26}. Moreover, we compute the spherical local relative character for this case.

We recall the notion of BZSV quadruples $\Delta=(G,H,\iota,\rho_H)$ in \cite{BZSV24} and further explicated in \cite{MWZ26b,TWZ26}, where $H$ is a closed subgroup of $G$, $\iota$ is a homomorphism from $\text{SL}_2$ into $G$ whose image commutes with $H$ and $\rho_H$ is a symplectic representation of $H$. We may consider $\iota$ also as a nilpotent orbit in $\mathfrak{g}=\text{Lie}(G)$ containing the image $d\iota\left(\begin{matrix}
0 & 1 \\
0 & 0
\end{matrix}\right)$. For this particular spherical variety $X$, the corresponding BZSV quadruple is $(F_4,\text{Spin}_9,0,0)$. As is also listed in \cite[Table 10]{TWZ26}, its BZSV dual quadruple is $\hat{\Delta}=(\hat{G},\hat{H}^\prime,\hat{\iota}^\prime,\rho_{\hat{H}^\prime})=(F_4,\text{SL}_2,C_3,0)$, noting that here we use the Bala-Carter label to represent the nilpotent orbit $\hat{\iota}^\prime$ in $\hat{\Delta}$. We denote by $\hat{\iota}^\prime_{BV}$ the Barbasch-Vogan dual nilpotent orbit of $\hat{\iota}^\prime$ in $F_4$. By \cite[Figure 2]{Ach03}, we have that the Bala-Carter label of $\hat{\iota}^\prime_{BV}$ is $A_2$.

Let $P_{BV}=M_{BV}N_{BV}$ be the parabolic subgroup of $G$ corresponding to $\hat{\iota}^\prime_{BV}$. Then $M_{BV}\cong \text{GSp}_6$. Let $\psi$ be a nontrivial additive character of $F$. This induces a character $\psi_{BV}$ of $N_{BV}(F)$. Recall the minimal representation $\Pi_\psi$ of $E_7$ in this case. By \cite[Corollary 3.7]{KS23}, we have the following isomorphism of $\text{PGL}_2(F)$-modules
$$
q:(\Pi_\psi^\infty)_{N_{BV},\psi_{BV}}\cong C^\infty(\Omega)_{N_{BV},\psi_{BV}}\cong C^\infty_c(\text{PGL}_2(F)).
$$
We consider the following commutative diagram.
\begin{equation}\label{BVdiagram}
\begin{tikzcd}
                                                                           & \Pi_\psi^\infty \arrow[ld, "q"'] \arrow[rd, "\theta_\sigma"] &                                                                            \\
C^\infty_c(\text{PGL}_2(F)) \arrow[d, "J"]                                 &                                                              & \sigma\otimes \pi:=\sigma\otimes \theta_\psi(\sigma) \arrow[ldd, "{l_{\pi,BV}}"] \\
{C^\infty_c(N,\psi\backslash \text{PGL}_2(F))} \arrow[rd, "\alpha_\sigma"] &                                                              &                                                                            \\
                                                                           & \sigma                                                       &                                                                           
\end{tikzcd}
\end{equation}
where $J$ is the map $f\mapsto \int_{N(F)}f(ng)\psi(n)^{-1}dn$, and $\pi=\theta_\psi(\sigma)$ is the theta lift of $\sigma$ and most importantly, $l_{\pi,BV}\in \text{Hom}_{N_{BV}}(\pi,\psi_{BV})$ is defined so that the above diagram commutes. The linear functional $l_{\pi,BV}$ is the desired local relative character of the degenerate Whittaker period $F_4/N_{BV}$, which resolves one of the key issues mentioned in Remark 1.7 of \cite{MWZ26}.

In the remainder of this section, we compute the spherical local relative character of the period $F_4/N_{BV}$. We adopt a similar set-up as in Section \ref{sec5}. Let $\Phi_0$ be the spherical vector of $\Pi_\psi$ as in Section \ref{sec5}. Assume that $\sigma$ is $\text{PGL}_2(\mathcal{O}_F)$-unramified. Then $\pi=\theta_\psi(\sigma)$ is $F_4(\mathcal{O}_F)$-unramified. We fix a $\text{PGL}_2(\mathcal{O}_F)$-spherical vector $v_0\in \sigma$ and $F_4(\mathcal{O}_F)$-spherical vector $w_0\in \pi$ such that $\langle v_0,v_0 \rangle_\sigma=\langle w_0,w_0 \rangle_\pi=1$. Then $\theta_\sigma(\Phi_0)=c_0v_0\otimes w_0$ for some $c_0\in \mathbb{C}$. From the diagram (\ref{BVdiagram}), we have
$$
\langle (\alpha_\sigma\circ J)(q(\Phi_0)),v_0 \rangle_\sigma = \langle l_{\pi,BV}(\theta_\sigma(\Phi_0)),v_0 \rangle_\sigma.
$$
Moreover, the left hand side is equal to
$$
\langle J\circ q(\Phi_0),\beta_{\sigma}(v_0) \rangle_{(N,\psi)\backslash \text{PGL}_2}=l_\sigma(v_0)\cdot \text{Vol}(\text{PGL}_2(\mathcal{O}_F)),
$$
as $(J\circ q)(\Phi_0)$ is the basic function of $(N,\psi)\backslash \text{PGL}_2$. Meanwhile, the right hand side is equal to $c_\sigma\cdot l_{\pi,BV}(w_0)$. By combining the above results with the computation in Proposition \ref{prop:unramified}, we obtain
$$
\left|l_{\pi,BV}(w_0)\right|^2=\frac{|l_\sigma(v_0)|^2\cdot \text{Vol}(\text{PGL}_2(\mathcal{O}_F))^2}{|c_\sigma|^2}=\frac{\zeta_F(4)\cdot\zeta_F(8)}{L(1,\sigma,\text{Ad})\cdot L(\frac{11}{2},\sigma,\text{std})\cdot L(\frac{5}{2},\sigma,\text{std})}.
$$
Therefore,
$$
\left|l_{\pi,BV}(w_0)\right|\cdot \left|l_{\pi}(w_0)\right|=\frac{1}{L(1,\sigma,\text{Ad})\cdot \text{Vol}(X(\mathcal{O}_F))}
$$
as desired. We summarize the above computation in the following theorem.
\begin{thm}Given $\Phi_0$ and $w_0$ defined as above, we have 
$$
\left|l_{\pi,BV}(w_0)\right|^2=\frac{\zeta_F(4)\cdot\zeta_F(8)}{L(1,\sigma,\text{Ad})\cdot L(\frac{11}{2},\sigma,\text{std})\cdot L(\frac{5}{2},\sigma,\text{std})}
$$
and
$$
\left|l_{\pi,BV}(w_0)\right|\cdot \left|l_{\pi}(w_0)\right|=\frac{1}{L(1,\sigma,\text{Ad})\cdot \text{Vol}(X(\mathcal{O}_F))},
$$
for any unramified representation $\sigma$, where $\pi=\theta_\psi(\sigma)$ and $l_{\pi,BV}$ is the local degenerate Whittaker functional defined in (\ref{BVdiagram}).
\end{thm}

\end{document}